\documentclass{birkmult}
 \newtheorem{thm}{Theorem}[section]
 \newtheorem{cor}[thm]{Corollary}
 
 \newtheorem{prop}[thm]{Proposition}
 \theoremstyle{definition}
 
 \theoremstyle{remark}
 \newtheorem{rem}[thm]{Remark}
 
 \numberwithin{equation}{section}

\begin{document}
%
%
%
%
%
%
%
%
%
\title[On the {Bessmertny\u{\i}} Class]
 {On the {Bessmertny\u{\i}} Class of Homogeneous\\ Positive Holomorphic Functions\\ on a Product of Matrix Halfplanes}
\author[Dmitry S. Kalyuzhny\u{\i}-Verbovetzki\u{\i}]{Dmitry S. Kalyuzhny\u{\i}-Verbovetzki\u{\i}}

\address{%
Department of Mathematics\\
Ben-Gurion University of the Negev\\
P.O. Box 653\\
Beer-Sheva 84105\\
Israel}

\email{dmitryk@math.bgu.ac.il}
\thanks{The   author was supported by the Center for Advanced Studies in Mathematics, Ben-Gurion University of the Negev.}



\date{}

\begin{abstract}
We generalize our earlier results from \cite{K} on the Bessmertny\u{\i} class of operator-valued functions holomorphic in the open right poly-halfplane which admit representation as a Schur complement of a block of a linear homogeneous operator-valued function with positive semidefinite operator coefficients, to the case of a product of open right matrix halfplanes. Several equivalent characterizations of this generalized Bessmertny\u{\i} class are presented. In particular, its intimate connection with the Agler--Schur class of holomorphic contractive operator-valued functions on the product of matrix unit disks is established.
\end{abstract}

\maketitle

\section{Introduction}\label{sec:intro}
In the PhD Thesis of M. F. Bessmertrny\u{\i} \cite{Bess} (the translation into English of some of its parts can be found in \cite{Bess1,Bess2,Bess3}) the class of rational $n\times n$ matrix-valued functions of $N$ complex variables $z=(z_1,\ldots,z_N)\in\mathbb{C}^N$, representable in the form
\begin{equation}\label{lrr}
f(z)=a(z)-b(z)d(z)^{-1}c(z),
\end{equation}
where a linear $(n+p)\times (n+p)$ matrix-valued function
\begin{equation}\label{pencil}
A(z)=z_1A_1+\cdots +z_NA_N=\left[\begin{array}{cc}
a(z) & b(z)\\
c(z) & d(z)
\end{array}\right]
\end{equation}
has positive semidefinite matrix coefficients $A_j,\ j=1,\ldots,N$, with real entries, was considered. Such a representation \eqref{lrr}--\eqref{pencil}, was called in the thesis a \emph{long resolvent representation}. The motivation of its consideration comes from  the electrical engineering. Bessmertrny\u{\i} has shown that this class is precisely the class of all characteristic functions of passive $2n$-poles, where the impedances of elements of an electrical circuit are considered as independent variables.

In \cite{K} a more general class $\mathcal{B}_N(\mathcal{U})$ of holomorphic functions on the \emph{open right poly-halfplane} $\Pi^N:=\{ z\in\mathbb{C}^N:\ {\rm Re}\, z_k> 0,\ k=1,\ldots,N\}$, with values in the $C^*$-algebra $L(\mathcal{U})$ of bounded linear operators on a Hilbert space $\mathcal{U}$,  which admit a representation \eqref{lrr} with a linear pencil $A(z)$ as in \eqref{pencil}, however consisting of operators from $L(\mathcal{U\oplus H})$ where $\mathcal{H}$ is an auxiliary Hilbert space, such that $A_j\geq 0,\ j=1,\ldots,N$, was introduced. Here the Hilbert spaces are supposed to be complex. This class $\mathcal{B}_N(\mathcal{U})$ was called the \emph{Bessmertrny\u{\i} class}. Any function $f\in\mathcal{B}_N(\mathcal{U})$ is homogeneous of degree one and takes operator values with positive semidefinite real parts. Moreover, $f$ can be uniquely extended to a holomorphic and homogeneous of degree one function on the domain \begin{equation}\label{omega}
\Omega_N:=\bigcup_{\lambda\in\mathbb{T}}(\lambda\Pi)^N\subset\mathbb{C}^N,
\end{equation}
so that \eqref{lrr} holds true for $z\in\Omega_N$, as well as the homogeneity relation
\begin{equation}\label{hom}
f(\lambda z)=\lambda f(z),\quad \lambda\in\mathbb{C}\setminus\{ 0\},\ z\in\Omega_N,
\end{equation}
and the symmetry relation
\begin{equation}\label{sym}
f(\bar{z})=f(z)^*,\quad z\in\Omega_N
\end{equation}
(here $\lambda z=(\lambda z_1,\ldots,\lambda z_N)$ and
$\bar{z}=(\bar{z}_1,\ldots,\bar{z}_N))$. In \cite{K} several
equivalent characterizations of the Bessmertny\u{\i} class have
been established: in terms of certain positive semidefinite
kernels on $\Omega_N\times\Omega_N$, in terms of functional
calculus of $N$-tuples of commuting bounded strictly accretive
operators on a common Hilbert space, and in terms of the double
Cayley transform. Let us briefly recall the last one. The
\emph{double Cayley transform} (over the variables and over the
values), $\mathcal{F}=\mathcal{C}(f)$, of a function
$f\in\mathcal{B}_N(\mathcal{U})$ is defined for $w$ in the open
unit polydisk $\mathbb{D}^N:=\{ w\in\mathbb{C}^N: \, |w_k|<1,\,
k=1,\ldots,N\}$ as
\begin{equation}\label{double-C}
\mathcal{F}(w)= \left(
f\left(\frac{1+w_1}{1-w_1},\ldots,\frac{1+w_N}{1-w_N}\right)
-I_\mathcal{U}\right)\left(
f\left(\frac{1+w_1}{1-w_1},\ldots,\frac{1+w_N}{1-w_N}\right)
+I_\mathcal{U}\right)^{-1}.
\end{equation}
For any $f\in\mathcal{B}_N(\mathcal{U})$, its double Cayley transform $\mathcal{F}=\mathcal{C}(f)$ belongs to the \emph{Agler--Schur class} $\mathcal{AS}_N(\mathcal{U})$, i.e., $\mathcal{F}$ is holomorphic on $\mathbb{D}^N$ and $\|\mathcal{F}(\mathbf{T})\|\leq 1$ for every $N$-tuple $\mathbf{T}=(T_1,\ldots,T_N)$ of commuting strict contractions on a common Hilbert space (see details on this class in \cite{Ag}). Moreover, there exist Hilbert spaces
        $\mathcal{X},\mathcal{X}_1,\ldots,\mathcal{X}_N$ such
        that $\mathcal{X}=\bigoplus_{k=1}^N\mathcal{X}_k$, and an
        \emph{Agler representation}
\begin{equation}\label{tf}
        \mathcal{F}(w)=D+CP(w)(I_\mathcal{X}-AP(w))^{-1}B,\quad w\in\mathbb{D}^N,
\end{equation}
    where $P(w)=\sum_{k=1}^Nw_kP_{\mathcal{X}_k}$, with orthogonal projections $P_{\mathcal{X}_k}$ onto $\mathcal{X}_k$, and
     $$\left[\begin{array}{cc}
       A & B \\
       C & D
     \end{array}\right]=:U=U^{-1}=U^*\in
     L(\mathcal{X}\oplus\mathcal{U}).$$
Conversely, any function $\mathcal{F}\in\mathcal{AS}_N(\mathcal{U})$ satisfying the latter condition can be represented as the double Cayley transform, $\mathcal{F}=\mathcal{C}(f)$, of some function
$f\in\mathcal{B}_N(\mathcal{U})$.

Let us recollect that matrices $A_j,\ j=1,\ldots,N$, in original Bessmertny\u{\i}'s definition had real entries, thus functions from his class took matrix values whose all entries were real at real points $z\in\mathbb{R}^N$. In \cite{K} we have considered also a ``real" version of the (generalized) Bessmertny\u{\i} class. Namely, we have  defined the real structure on a Hilbert space $\mathcal{U}$ by means of  an \emph{anti-unitary involution} (a counterpart of the complex conjugation), i.e., an operator  $\iota=\iota_\mathcal{U}:\mathcal{U}\rightarrow\mathcal{U}$ such that
\begin{eqnarray}
\iota^2 & = & I_{\mathcal U},\label{invol}\\
\left\langle \iota u_1,\iota u_2\right\rangle &=& \left\langle u_2,u_1\right\rangle,\quad u_1,u_2\in\mathcal{U}.\label{antiun}
\end{eqnarray}
Such an operator $\iota$ is \emph{anti-linear}, i.e.,
$$\iota(\alpha u_1+\beta u_2)=\bar{\alpha}u_1+\bar{\beta}u_2,\quad \alpha,\beta\in\mathbb{C},\ u_1,u_2\in\mathcal{U}.$$
An operator $A$ from $L(\mathcal{U,Y})$, the Banach space of all bounded linear operators from a Hilbert space $\mathcal{U}$ to a Hilbert space $\mathcal{Y}$, is called \emph{$(\iota_{\mathcal U},\iota_{\mathcal Y})$-real} for anti-unitary involutions $\iota_{\mathcal U}$ and $\iota_{\mathcal Y}$ if
\begin{equation}\label{real-op}
\iota_{\mathcal Y}A=A\iota_{\mathcal U}.
\end{equation}
Such operators $A$ are a counterpart of matrices with real entries. Finally, a function $f$ on a set $\Omega\subset\mathbb{C}^N$ such that $z\in\Omega\Leftrightarrow \bar{z}\in\Omega$, which takes values from $L(\mathcal{U,Y})$  is called \emph{$(\iota_{\mathcal U},\iota_{\mathcal Y})$-real} if
\begin{equation}\label{real-func}
f^\sharp (z):=\iota_{\mathcal Y}f(\bar{z})\iota_{\mathcal U}=f(z),\quad z\in\Omega.
\end{equation}
If $\mathcal{U}=\mathcal{Y}$ and $\iota_{\mathcal U}=\iota_{\mathcal Y}=\iota$ then such a function is called \emph{$\iota$-real}. We have defined the ``\emph{$\iota$-real" Bessmertny\u{\i} class} $\iota\mathbb{R}\mathcal{B}_N(\mathcal{U})$ as the subclass of all $\iota$-real functions from $\mathcal{B}_N(\mathcal{U})$. The latter subclass is a counterpart of the original class considered by Bessmertny\u{\i}.
In \cite{K} we have obtained different characterizations for $\iota\mathbb{R}\mathcal{B}_N(\mathcal{U})$, too.

In the present paper we introduce and investigate analogous classes of functions (either for the ``complex" and ``real" cases) on more general domains. First, we define a \emph{product of matrix halfplanes} as
\begin{eqnarray}
\Pi^{n_1\times n_1}\times\cdots\times\Pi^{n_N\times n_N}:=\{ Z=(Z_1,\ldots,Z_N):\nonumber \\
 Z_k\in\mathbb{C}^{n_k\times n_k},\ Z_k+Z_k^*>0,\quad k=1,\ldots,N\}\label{product}
\end{eqnarray}
 which serves as a generalization of the open right poly-halfplane $\Pi^N$. Then we define a conterpart of the domain $\Omega_N$ as
 \begin{equation}\label{omega-g}
\Omega_{n_1,\ldots,n_N}:=\bigcup_{\lambda\in\mathbb{T}}\left(\lambda\Pi^{n_1\times n_1}\times\cdots\times\lambda\Pi^{n_N\times n_N}\right),
\end{equation}
and define the corresponding Bessmertny\u{\i} classes of functions on the domain $\Omega_{n_1,\ldots,n_N}$. Consideration of such classes can be also motivated by problems of the theory of electrical networks since there are situations where ``matrix impedances" are considered as matrix variables (see, e.g., \cite{Kr}). On the other hand, mathematical tools for such an investigation have recently appeared. Since in \cite{K} the close relation of the Bessmertny\u{\i} classes $\mathcal{B}_N(\mathcal{U})$ and
$\iota\mathbb{R}\mathcal{B}_N(\mathcal{U})$ to the Agler--Schur class $\mathcal{AS}_N(\mathcal{U})$ has been established, this has made possible the use of properties of the latter class as a tool for investigation Bessmertny\u{\i}'s classes. In the same manner we make use of the recent works of C.-G.~Ambrozie and D.~Timotin \cite{AT}, J.~A.~Ball and V.~Bolotnikov \cite{BB} on the Agler--Schur class of function on so-called polynomially defined domains for the investigation of the Bessmertny\u{\i}'s classes of functions on $\Omega_{n_1,\ldots,n_N}$. A counterpart of the class
$\mathcal{B}_N(\mathcal{U})$ is introduced in Section~\ref{sec:b}, where also a useful decomposition for functions from this class is obtained. In Section~\ref{sec:b-calc} the relationship between the Bessmertny\u{\i} class on $\Omega_{n_1,\ldots,n_N}$ and the corresponding Agler--Schur class on a product of matrix disks is established. This allows us to give a characterization of the (generalized) Bessmertny\u{\i} class in terms of functional calculus for collections of operators. In Section~\ref{sec:image} we describe the image of this class under the double Cayley transform. Finally, a counterpart of the ``real" Bessmertny\u{\i} class
$\iota\mathbb{R}\mathcal{B}_N(\mathcal{U})$ is studied in Section~\ref{sec:real-b}.

\section{The Bessmertny\u{\i} class for a matrix domain}\label{sec:b}
Let us define the class $\mathcal{B}_{n_1,\ldots,n_N}(\mathcal{U})$ of all $L(\mathcal{U})$-valued functions $f$ holomorphic on the domain $\Omega_{n_1,\ldots,n_N}$ defined in \eqref{omega-g} (see also \eqref{product}) which are representable as
\begin{equation}\label{sc}
f(Z)=a(Z)-b(Z)d(Z)^{-1}c(Z)
\end{equation}
for $Z\in \Omega_{n_1,\ldots,n_N},$ where
\begin{equation}\label{lp}
A(Z)=G_1^*(Z_1\otimes I_{\mathcal{M}_1})G_1+\cdots +G_N^*(Z_N\otimes I_{\mathcal{M}_N})G_N=
\left[\begin{array}{cc}
a(Z) & b(Z)\\
c(Z) & d(Z)
\end{array}\right]\in L(\mathcal{U\oplus H})
\end{equation}
for some Hilbert spaces $\mathcal{M}_1,\ldots,\mathcal{M}_N,\mathcal{H}$ and operators $G_k\in L(\mathcal{U\oplus H},\mathbb{C}^{n_k}\otimes\mathcal{M}_k),\ k=1,\ldots,N$.
\begin{rem}\label{rem:dom}
If a function $f$ is holomorphic on $\Pi^{n_1\times n_1}\times\cdots\times\Pi^{n_N\times n_N}$ and has a representation \eqref{sc}--\eqref{lp} there, then $f$ can be extended to
$\Omega_{n_1,\ldots,n_N}$ by homogeneity of degree one, and this extension is,
clearly, holomorphic and admits a representation \ref{sc} in
$\Omega_{n_1,\ldots,n_N}$. That is why we define the class
$\mathcal{B}_{n_1,\ldots,n_N}(\mathcal{U})$ straight away as a class of functions
on $\Omega_{n_1,\ldots,n_N}$. Keeping in mind the possibility and uniqueness of
such extension, we will write sometimes
$f\in\mathcal{B}_{n_1,\ldots,n_N}(\mathcal{U})$ for functions defined originally
on $\Pi^{n_1\times n_1}\times\cdots\times\Pi^{n_N\times n_N}$.
\end{rem}
\begin{thm}\label{thm:b-decomp}
Let $f$ be an $L\mathcal{(U)}$-valued function holomorphic on
$\Pi^{n_1\times n_1}\times\cdots\times\Pi^{n_N\times n_N}$. Then $f\in\mathcal{B}_{n_1,\ldots,n_N}(\mathcal{U})$ if and only if
there exist holomorphic functions  $\varphi_k(Z)$ on $\Pi^{n_1\times n_1}\times\cdots\times\Pi^{n_N\times n_N}$ with values in $L(\mathcal{U},\mathbb{C}^{n_k}\otimes\mathcal{M}_k),\
k=1,\ldots,N$, such that
\begin{equation}\label{b-decomp}
    f(Z)=\sum\limits_{k=1}^N\varphi_k(\Lambda)^*
    (Z_k\otimes I_{\mathcal{M}_k})\varphi_k(Z),\quad Z,\Lambda \in\Pi^{n_1\times n_1}\times\cdots\times\Pi^{n_N\times n_N}
\end{equation}
holds. In this case the functions $\varphi_k(Z)$
can be uniquely extended to the holomorphic functions on
$\Omega_{n_1,\ldots,n_N}$ (we use the same notation for the
extended functions) which are homogeneous of degree zero, i.e., for
every $\lambda\in\mathbb{C}\backslash\{ 0\}$,
\begin{equation}\label{homzero}
    \varphi_k(\lambda Z)=\varphi_k(Z),\quad
    Z\in\Omega_{n_1,\ldots,n_N},
\end{equation}
 and identity \eqref{b-decomp} is extended to all of
$Z,\Lambda\in\Omega_{n_1,\ldots,n_N}$.
\end{thm}
\begin{proof}
\textbf{Necessity.}
Let $f\in\mathcal{B}_{n_1,\ldots,n_N}(\mathcal{U})$. Then \eqref{sc} holds for $Z\in\Omega_{n_1,\ldots,n_N}$,
 some Hilbert spaces $\mathcal{H},\mathcal{M}_1,\ldots,\mathcal{M}_N$ and a linear pencil of
operators  \eqref{lp}. Define
$$\psi(Z):=\left[\begin{array}{c}
  I_\mathcal{U} \\
  -d(Z)^{-1}c(Z)
\end{array}\right]\in L\mathcal{(U,U\oplus H)},\quad Z\in\Omega_{n_1,\ldots,n_N}.$$
Then for all $Z,\Lambda\in\Omega_{n_1,\ldots,n_N}$ one has
\begin{eqnarray*}
f(Z) &=& a(Z)-b(Z)d(Z)^{-1}c(Z)\\
     &=& \left[\begin{array}{cc}
       I_\mathcal{U} & -c(\Lambda)^*d(\Lambda)^{-*}
     \end{array}\right]\left[\begin{array}{c}
       a(Z)-b(Z)d(Z)^{-1}c(Z) \\
       0
     \end{array}\right]\\
     &=& \left[\begin{array}{cc}
       I_\mathcal{U} & -c(\Lambda)^*d(\Lambda)^{-*}
     \end{array}\right]\left[\begin{array}{cc}
       a(Z) & b(Z) \\
       c(Z) & d(Z)
     \end{array}\right]\left[\begin{array}{c}
       I_\mathcal{U} \\
       -d(Z)^{-1}c(Z)
     \end{array}\right]\\
     &=& \psi(\Lambda)^*A(Z)\psi(Z).
\end{eqnarray*}
Set $\varphi_k(Z):=G_k\psi(Z),\ k=1,\ldots,N$. Clearly, the
functions $\varphi_k(Z),\ k=1,\ldots,N$, are holomorphic on
$\Omega_{n_1,\ldots,n_N}$ and satisfy \eqref{homzero}. Rewriting
the equality
\begin{equation}\label{f-id}
f(Z)=\psi(\Lambda)^*A(Z)\psi(Z),\quad Z\in\Omega_{n_1,\ldots,n_N},
\end{equation}
in the form
\begin{equation}\label{b-decomp-ext}
    f(Z)=\sum_{k=1}^N\varphi_k(\Lambda)^*(Z_k\otimes I_{\mathcal{M}_k})\varphi_k(Z),\quad Z,\Lambda\in\Omega_{n_1,\ldots,n_N},
\end{equation}
we obtain, in particular, \eqref{b-decomp}.

\textbf{Sufficiency.}
Let $f$ be an $L\mathcal{(U)}$-valued function holomorphic on
$\Pi^{n_1\times n_1}\times\cdots\times\Pi^{n_N\times n_N}$ and representable there in the form \eqref{b-decomp} with
some holomorphic functions $\varphi_k(Z)$ taking values in $L(\mathcal{U},\mathbb{C}^{n_k}\otimes\mathcal{M}_k),\ k=1,\ldots,N$. Set
$$\mathcal{N}:=\bigoplus_{k=1}^N(\mathbb{C}^{n_k}\otimes\mathcal{M}_k),\
P_k:=P_{\mathcal{M}_k},\quad k=1,\ldots,N,$$
$$\varphi(Z):=\mbox{col}\left(\begin{array}{ccc}
  \varphi_1(Z) & \ldots & \varphi_N(Z)
\end{array}\right)\in L\mathcal{(U,N)},$$
$$E:=\left(I_{n_1},\ldots,I_{n_N}\right)\in\Pi^{n_1\times n_1}\times\cdots\times\Pi^{n_N\times n_N},$$
where $I_n$ denotes the identity $n\times n$ matrix.
From \eqref{b-decomp} we get
\begin{equation}\label{e1}
    f(E)=\sum_{k=1}^N\varphi_k(\Lambda)^*\varphi_k(E),\quad
    \Lambda\in\Pi^{n_1\times n_1}\times\cdots\times\Pi^{n_N\times n_N}.
\end{equation}
In particular,
\begin{equation}\label{e2}
    f(E)=\sum_{k=1}^N\varphi_k(E)^*\varphi_k(E).
\end{equation}
By subtracting \eqref{e2} from \eqref{e1} we get
$$\sum_{k=1}^N[\varphi_k(\Lambda)-\varphi_k(E)]^*\varphi_k(E)=0,\quad
    \Lambda\in\Pi^{n_1\times n_1}\times\cdots\times\Pi^{n_N\times n_N},$$
i.e., the following orthogonality relation holds:
$$\mathcal{H}:=\mbox{clos span}_{\Lambda\in\Pi^{n_1\times n_1}\times\cdots\times\Pi^{n_N\times n_N}}\{
[\varphi(\Lambda)-\varphi(E)]\mathcal{U}\}
\perp\mbox{clos}\{\varphi(E)\mathcal{U}\} =:\mathcal{X}.$$ For any
$\Lambda\in\Pi^{n_1\times n_1}\times\cdots\times\Pi^{n_N\times n_N}$ and $u\in\mathcal{U}$ one can represent now
$\varphi(\Lambda)u$ as $$\mbox{col}\left[\begin{array}{cc}
  \varphi(E) & \varphi(\Lambda)-\varphi(E)
\end{array}\right]u\in\mathcal{X\oplus H}.$$ On the other hand, for
any $u\in\mathcal{U},\ \Lambda\in\Pi^{n_1\times n_1}\times\cdots\times\Pi^{n_N\times n_N}$ one has
\begin{eqnarray*}
\varphi (E)u & \in &\mbox{clos span}_{\Lambda\in\Pi^{n_1\times n_1}\times\cdots\times\Pi^{n_N\times n_N}}\{ \varphi(\Lambda )\mathcal{U}\}, \\
(\varphi(\Lambda )-\varphi(E))u &\in &\mbox{clos span}_{\Lambda\in\Pi^{n_1\times n_1}\times\cdots\times\Pi^{n_N\times n_N}}\{ \varphi(\Lambda )\mathcal{U}\}.
\end{eqnarray*}
 Thus, $\mbox{clos
span}_{\Lambda\in\Pi^{n_1\times n_1}\times\cdots\times\Pi^{n_N\times n_N}}\{\varphi(\Lambda )\mathcal{U}\}=\mathcal{X\oplus
H}$. Let $\kappa:\mathcal{X\oplus H}\rightarrow\mathcal{N}$ be
the natural embedding defined by
\begin{equation}\label{kappa}
   \kappa:\left[\begin{array}{c}
  \varphi(E)u \\
(\varphi(\Lambda )-\varphi(E))u
\end{array}\right]\longmapsto\varphi(\Lambda )u=\left[\begin{array}{c}
  \varphi_1(\Lambda )u \\
  \vdots \\
  \varphi_N(\Lambda )u
\end{array}\right]
\end{equation}
and extended to the whole $\mathcal{X\oplus H}$ by linearity and
continuity. Set
$$G_k:=(I_{n_k}\otimes P_k)\kappa\left[\begin{array}{cc}
  \varphi(E) & 0 \\
  0 & I_\mathcal{H}
\end{array}\right]\in L(\mathcal{U\oplus H},\mathbb{C}^{n_k}\otimes\mathcal{M}_k),\quad k=1,\ldots,N,$$
$$\psi(\Lambda):=\left[\begin{array}{c}
  I_\mathcal{U}\\
\varphi(\Lambda )-\varphi(E)
\end{array}\right]\in L\mathcal{(U,U\oplus H)},\quad
\Lambda\in\Pi^{n_1\times n_1}\times\cdots\times\Pi^{n_N\times n_N}.$$ Then $$f(Z)=\psi(\Lambda)^*A(Z)\psi(Z),\quad
Z,\Lambda\in\Pi^{n_1\times n_1}\times\cdots\times\Pi^{n_N\times n_N},$$
where $A(Z)$ is defined by \eqref{lp}. Indeed,
\begin{eqnarray*}
\lefteqn{\psi(\Lambda)^*A(Z)\psi(Z)= \left[\begin{array}{c}
  I_\mathcal{U} \\
\varphi(\Lambda)-\varphi(E)
\end{array}\right]^*\left[\begin{array}{cc}
  \varphi(E) & 0 \\
  0 & I_\mathcal{H}
\end{array}\right]^*\kappa^*\left(\sum_{k=1}^NZ_k\otimes P_k\right)\kappa}\\
&\times
\left[\begin{array}{cc}
  \varphi(E) & 0 \\
  0 & I_\mathcal{H}
\end{array}\right]\left[\begin{array}{c}
  I_\mathcal{U} \\
\varphi(Z)-\varphi(E)
\end{array}\right]
=\sum\limits_{k=1}^N\varphi_k(\Lambda)^*(Z_k\otimes I_{\mathcal{M}_k})\varphi_k(Z)=f(Z).
\end{eqnarray*}
Now, with respect to the block partitioning of $A(Z)$ we have
\begin{eqnarray*}
A(Z)\psi(Z) &=& \left[\begin{array}{cc}
  a(Z) & b(Z) \\
  c(Z) & d(Z)
\end{array}\right]\left[\begin{array}{c}
  I_\mathcal{U} \\
\varphi(Z)-\varphi(E)
\end{array}\right]\\
&=& \left[\begin{array}{c}
  a(Z)+b(Z)(\varphi(Z)-\varphi(E)) \\
  c(Z)+d(Z)(\varphi(Z)-\varphi(E))
\end{array}\right]=:\left[\begin{array}{c}
  f_1(Z) \\
  f_2(Z)
\end{array}\right].
\end{eqnarray*}
Since for $Z,\Lambda\in\Pi^{n_1\times n_1}\times\cdots\times\Pi^{n_N\times n_N}$ one has
$$\psi(\Lambda)^*A(Z)\psi(Z)=\left[\begin{array}{cc}
  I_\mathcal{U} & \varphi(\Lambda)^*-\varphi(E)^*
\end{array}\right]\left[\begin{array}{c}
  f_1(Z) \\
  f_2(Z)
\end{array}\right]=f(Z),$$
by setting $\Lambda :=E$ in this equality we get $$f_1(Z)=f(Z),\quad
Z\in\Pi^{n_1\times n_1}\times\cdots\times\Pi^{n_N\times n_N}.$$ Therefore, for every $Z,\Lambda\in\Pi^{n_1\times n_1}\times\cdots\times\Pi^{n_N\times n_N}$
we get $[\varphi(\Lambda)-\varphi(E)]^*f_2(Z)=0$. This implies that
for every $Z\in\Pi^{n_1\times n_1}\times\cdots\times\Pi^{n_N\times n_N}$ and $u\in\mathcal{U}$ one has
$f_2(Z)u\perp\mathcal{H}$. But $f_2(Z)u\in\mathcal{H}$. Therefore,
$f_2(Z)u=0$, and $f_2(Z)\equiv 0$, i.e.,
\begin{equation}\label{id-vanish}
c(Z)+d(Z)[\varphi(Z)-\varphi(E)]\equiv 0.
\end{equation}
Since for every $Z\in\Pi^{n_1\times n_1}\times\cdots\times\Pi^{n_N\times n_N}$ the operator
$P(Z):=\sum_{k=1}^NZ_k\otimes P_k$ has positive definite real part, i.e.,
$P(Z)+P(Z)^*\geq\alpha_ZI_\mathcal{N}>0$ for some scalar $\alpha_Z>0$, the
operator $d(Z)=P_\mathcal{H}\kappa^*P(Z)\kappa |\mathcal{H}$ has
positive definite real part, too. Therefore, $d(Z)$ is boundedly
invertible for all $Z\in\Pi^{n_1\times n_1}\times\cdots\times\Pi^{n_N\times n_N}$. From \eqref{id-vanish} we get
$\varphi(Z)-\varphi(E)=-d(Z)^{-1}c(Z),\ Z\in\Pi^{n_1\times n_1}\times\cdots\times\Pi^{n_N\times n_N}$, and
$$f(Z)=f_1(Z)=a(Z)-b(Z)d(Z)^{-1}c(Z),\quad Z\in\Pi^{n_1\times n_1}\times\cdots\times\Pi^{n_N\times n_N}.$$
Taking into account Remark~\ref{rem:dom}, we get
$f\in\mathcal{B}_{n_1,\ldots,n_N}(\mathcal{U})$.

Functions $\varphi(Z)-\varphi(E)=-d(Z)^{-1}c(Z)$ and, hence, $\psi(Z)$
are well-defined, holomorphic and homogeneous of degree zero
 on $\Omega_{n_1,\ldots,n_N}$, thus \eqref{b-decomp-ext}
holds.

The proof is complete.
\end{proof}

\section{The class $\mathcal{B}_{n_1,\ldots,n_N}(\mathcal{U})$ and functional calculus}\label{sec:b-calc}
Let us observe now that \eqref{b-decomp} is equivalent to the couple of identities
\begin{eqnarray}\label{b-decomp-1}
    f(Z)+f(\Lambda)^* &=& \sum\limits_{k=1}^N\varphi_k(\Lambda)^*
    ((Z_k+\Lambda_k^*)\otimes I_{\mathcal{M}_k})\varphi_k(Z),\\
    \label{b-decomp-2}
    f(Z)-f(\Lambda)^* &=& \sum\limits_{k=1}^N\varphi_k(\Lambda)^*
    ((Z_k-\Lambda_k^*)\otimes I_{\mathcal{M}_k})\varphi_k(Z)
\end{eqnarray}
valid for all $Z,\Lambda \in\Pi^{n_1\times n_1}\times\cdots\times\Pi^{n_N\times n_N}$.
We will show that the double Cayley transform $\mathcal{F}=\mathcal{C}(f)$ applied to a function $f$ from the Bessmertny\u{\i} class $\mathcal{B}_{n_1,\ldots,n_N}(\mathcal{U})$ and defined as
\begin{eqnarray}\nonumber
\lefteqn{\mathcal{F}(W) = \left[
f((I_{n_1}+W_1)(I_{n_1}-W_1)^{-1},\ldots,(I_{n_N}+W_N)(I_{n_N}-W_N)^{-1})
-I_\mathcal{U}\right]}\\
&\times
\left[
f((I_{n_1}+W_1)(I_{n_1}-W_1)^{-1},\ldots,(I_{n_N}+W_N)(I_{n_N}-W_N)^{-1})
+I_\mathcal{U}\right]^{-1}\label{double-C-matr}
\end{eqnarray}
(compare with \eqref{double-C})
turns the first of these identities into an Agler-type identity which characterizes
the Agler--Schur class of holomorphic $L(\mathcal{U})$-valued functions on the \emph{product of open matrix unit disks}
\begin{eqnarray*}
\mathbb{D}^{n_1\times n_1}\times\cdots\times\mathbb{D}^{n_N\times n_N} &:=& \{ W=(W_1,\ldots,W_N)\in \mathbb{C}^{n_1\times n_1}\times\cdots\times\mathbb{C}^{n_N\times n_N}: \\
& & W_kW_k^*< I_{n_k},\quad k=1,\ldots,N\}.
\end{eqnarray*}
 The latter is a special case of the Agler--Schur class of holomorphic $L(\mathcal{U})$-valued functions on
a domain with matrix polynomial defining function, which was
studied in \cite{AT} and \cite{BB}. This allows us to obtain one
more characterization of
$\mathcal{B}_{n_1,\ldots,n_N}(\mathcal{U})$. Let $P(w),\
w\in\mathbb{C}^n$, be a polynomial $p\times q$ matrix valued
function, and $$\mathcal{D}_P:=\{ w\in\mathbb{C}^n:\ \| P(w)\|
<1\}$$ (here and in the sequel the norm of a $p\times q$ matrix
means its operator norm with respect to the standard Euclidean
metrics in $\mathbb{C}^p$ and $\mathbb{C}^q$). Let
$\mathcal{C}_{\mathcal{D}_P}$ denote the set of commutative
$n$-tuples $\mathbf{T}=(T_1,\ldots,T_n)$ of bounded linear
operators on a common Hilbert space $\mathcal{H}_\mathbf{T}$
subject to the condition $\| P(\mathbf{T})\| <1$. It was shown in
\cite{AT} that the \emph{Taylor joint spectrum}
$\sigma_T(\mathbf{T})$ (see \cite{T1,V1} and also \cite{Cu}) of
any $\mathbf{T}\in\mathcal{C}_{\mathcal{D}_P}$ is contained in
$\mathcal{D}_P$. Thus, for any function $S$ holomorphic on
$\mathcal{D}_P$ and any $\mathbf{T}\in\mathcal{C}_{\mathcal{D}_P}$
the operator $S(\mathbf{T})$ is well defined by the \emph{Taylor
functional calculus} (see \cite{T2,V2} and also \cite{Cu}). For
the domain $\mathcal{D}_P$, the \emph{Agler--Schur class}
$\mathcal{AS}_{\mathcal{D}_P}(\mathcal{E,E}_*)$ consists of all
holomorphic $L(\mathcal{E,E}_*)$-valued functions $\mathcal{F}$ on
$\mathcal{D}_P$ such that
\begin{equation}\label{vN}
\|\mathcal{F}(\mathbf{T})\|\leq 1,\quad \mathbf{T}\in\mathcal{C}_{\mathcal{D}_P}.
\end{equation}
Recall the following theorem from \cite{BB} (the case when $\mathcal{E}=\mathcal{E}_*=\mathbb{C}$ can be found in \cite{AT}), however in a slightly simplified form which will be sufficient for our purpose.
\begin{thm}\label{thm:BB}
Let $\mathcal{F}$ be an $L(\mathcal{E,E}_*)$-valued function holomorphic on $\mathcal{D}_P$. Then the following statements are equivalent:
\begin{description}
    \item[(i)] $\mathcal{F}\in\mathcal{AS}_{\mathcal{D}_P}(\mathcal{E,E}_*)$;
    \item[(ii)] there exist an auxiliary Hilbert space $\mathcal{M}$ and an $L(\mathbb{C}^p\otimes\mathcal{M,E}_*)$-valued  function $H^L$ holomorphic on $\mathcal{D}_P$ such that
    \begin{equation}\label{bb-id-l}
    I_{\mathcal{E}_*}-\mathcal{F}(w)\mathcal{F}(\omega)^*=H^L(w)\left((I_p-P(w)P(\omega)^*)\otimes I_\mathcal{M}\right)H^L(\omega)^*
    \end{equation}
    holds for all $w,\omega\in\mathcal{D}_P$;
    \item[(iii)] there exist an auxiliary Hilbert space $\mathcal{M}$ and an $L(\mathcal{E},\mathbb{C}^q\otimes\mathcal{M})$-valued function $H^R$  holomorphic on $\mathcal{D}_P$ such that
    \begin{equation}\label{bb-id-r}
I_{\mathcal{E}}-\mathcal{F}(\omega)^*\mathcal{F}(w)=H^R(\omega)^*\left((I_q-P(\omega)^*P(w))\otimes I_\mathcal{M}\right)H^R(w)
    \end{equation}
    holds for all $w,\omega\in\mathcal{D}_P$;
    \item[(iv)] there exist an auxiliary Hilbert space $\mathcal{M}$, an $L(\mathbb{C}^p\otimes\mathcal{M,E}_*)$-valued function $H^L$ and an $L(\mathcal{E},\mathbb{C}^q\otimes\mathcal{M})$-valued  function $H^R$, which are holomorphic on $\mathcal{D}_P$, such that
    \begin{eqnarray}
\left[\begin{array}{cc}
I_{\mathcal{E}}-\mathcal{F}(\omega')^*\mathcal{F}(w) & \mathcal{F}(\omega')^*-\mathcal{F}(\omega)^*\\
\mathcal{F}(w')-\mathcal{F}(w) & I_{\mathcal{E}_*}-\mathcal{F}(w')\mathcal{F}(\omega)^*
\end{array}\right]=\left[\begin{array}{cc}
H^R(\omega')^* & 0\\
0 & H^L(w')
\end{array}\right] \label{bb-id-lr}\\
\times\left(\left[\begin{array}{cc}
I_q-P(\omega')^*P(w) & P(\omega')^*-P(\omega)^* \\
P(w')-P(w) & I_p-P(w')P(\omega)^*
\end{array}\right]\otimes I_\mathcal{M}\right)\left[\begin{array}{cc}
H^R(w) & 0\\
0 & H^L(\omega)^*
\end{array}\right]\nonumber
    \end{eqnarray}
    holds for all $w,w',\omega,\omega'\in\mathcal{D}_P$;
    \item[(v)] there exists a Hilbert space $\mathcal{X}$ and a unitary operator
    \begin{equation}\label{bb-u}
    U=\left[\begin{array}{cc}
A & B\\
C & D
\end{array}\right]\in L((\mathbb{C}^p\otimes\mathcal{X})\oplus\mathcal{E},(\mathbb{C}^q\otimes\mathcal{X})\oplus\mathcal{E}_*)
\end{equation}
such that
\begin{equation}\label{bb-tf}
F(w)=D+C(P(w)\otimes I_\mathcal{X})\left(I_{\mathbb{C}^q\otimes\mathcal{X}}-A(P(w)\otimes I_\mathcal{X})\right)^{-1}B
\end{equation}
holds for all $w\in\mathcal{D}_P$.
\end{description}
\end{thm}
In \cite{BB} it was shown how to obtain from \eqref{bb-id-l} a unitary operator \eqref{bb-u} which gives the representation \eqref{bb-tf} for an arbitrary $\mathcal{F}\in\mathcal{AS}_{\mathcal{D}_P}(\mathcal{E,E}_*)$. We will show now how to get from \eqref{bb-id-lr}  a special unitary operator \eqref{bb-u} and representation \eqref{bb-tf} for an arbitrary $\mathcal{F}\in\mathcal{AS}_{\mathcal{D}_P}(\mathcal{E,E}_*)$. Let \eqref{bb-id-lr} hold for such $\mathcal{F}$, where a Hilbert space $\mathcal{M}$ and functions $H^L,H^R$ are such as in statement (iv) of Theorem~\ref{thm:BB}. Define the lineals
\begin{eqnarray*}
\mathcal{D}_0 & := & \mbox{span}\,\left\{\left[\begin{array}{c}
(P(w)\otimes I_\mathcal{M})H^R(w)\\
I_\mathcal{E}\end{array}\right] e,\left[\begin{array}{c}
H^L(\omega)^*\\
\mathcal{F}(\omega)^*
\end{array}\right] e_*:\right.\\
& & \left. w,\omega\in\mathcal{D}_P,\ e\in\mathcal{E},\ e_*\in\mathcal{E}_*\ \right\} \subset(\mathbb{C}^p\otimes\mathcal{M})\oplus\mathcal{E},  \\
\mathcal{R}_0 & := & \mbox{span}\,\left\{\left[\begin{array}{c}
H^R(w)\\
\mathcal{F}(w)\end{array}\right] e,\left[\begin{array}{c}
(P(\omega)^*\otimes I_\mathcal{M})H^L(\omega)^*\\
I_{\mathcal{E}_*}
\end{array}\right] e_*:\right.\\
& & \left. w,\omega\in\mathcal{D}_P,\ e\in\mathcal{E},\ e_*\in\mathcal{E}_*\ \right\} \subset(\mathbb{C}^q\otimes\mathcal{M})\oplus\mathcal{E}_*,
\end{eqnarray*}
and the operator $U_0:\mathcal{D}_0\rightarrow \mathcal{R}_0$ which acts on the generating vectors of $\mathcal{D}_0$ as
\begin{eqnarray*}
 \left[\begin{array}{c}
(P(w)\otimes I_\mathcal{M})H^R(w)\\
I_\mathcal{E}\end{array}\right] e \longmapsto\left[\begin{array}{c}
H^R(w)\\
\mathcal{F}(w)\end{array}\right] e,\quad w\in\mathcal{D}_P,\ e\in\mathcal{E},\\
 \left[\begin{array}{c}
H^L(\omega)^*\\
\mathcal{F}(\omega)^*
\end{array}\right] e_*\longmapsto\left[\begin{array}{c}
(P(\omega)^*\otimes I_\mathcal{M})H^L(\omega)^*\\
I_{\mathcal{E}_*}
\end{array}\right] e_*,\quad \omega\in\mathcal{D}_P,\  e_*\in\mathcal{E}_*.
\end{eqnarray*}
This operator is correctly defined. Moreover, $U_0$ maps $\mathcal{D}_0$ isometrically onto $\mathcal{R}_0$. Indeed,
\eqref{bb-id-lr} can be rewritten as
\begin{eqnarray*}
& \left[\begin{array}{cc}
H^R(\omega')^* & \mathcal{F}(\omega')^*\\
H^L(w')(P(w')\otimes I_\mathcal{M}) & I_{\mathcal{E}_*}
\end{array}\right]\left[\begin{array}{cc}
H^R(w) & (P(\omega)^*\otimes I_\mathcal{M})H^L(\omega)^*\\
\mathcal{F}(w) &   I_{\mathcal{E}_*}
\end{array}\right]\\
=& \left[\begin{array}{cc}
H^R(\omega')^*(P(\omega')^*\otimes I_\mathcal{M}) & I_{\mathcal{E}}\\
H^L(w') & F(w')
\end{array}\right]\left[\begin{array}{cc}
(P(w)\otimes I_\mathcal{M}) H^R(w) & H^L(\omega)^*\\
I_\mathcal{E} &   \mathcal{F}(\omega)^*
\end{array}\right],
    \end{eqnarray*}
which means that for
\begin{eqnarray*}
x &=& \left[\begin{array}{c}
(P(w)\otimes I_\mathcal{M})H^R(w)\\
I_\mathcal{E}\end{array}\right] e+\left[\begin{array}{c}
H^L(\omega)^*\\
\mathcal{F}(\omega)^*
\end{array}\right] e_*,\\
x'&=& \left[\begin{array}{c}
(P(w')\otimes I_\mathcal{M})H^R(w')\\
I_\mathcal{E}\end{array}\right] e'+\left[\begin{array}{c}
H^L(\omega')^*\\
\mathcal{F}(\omega')^*
\end{array}\right] e'_*,
\end{eqnarray*}
one has
$$\left\langle U_0x,U_0x'\right\rangle =\left\langle x,x'\right\rangle.$$
Clearly, $U_0$ can be uniquely extended to the unitary operator $\widetilde{U_0}:\ \mbox{clos}(\mathcal{D}_0)\to
\mbox{clos}(\mathcal{R}_0)$. In the case when
\begin{equation}\label{dim}
\dim\{ ((\mathbb{C}^p\otimes\mathcal{M})\oplus\mathcal{E})\ominus\mbox{clos}(\mathcal{D}_0)\} =
\dim\{ ((\mathbb{C}^q\otimes\mathcal{M})\oplus\mathcal{E}_*)\ominus\mbox{clos}(\mathcal{R}_0)\}
\end{equation}
there exists a (non-unique!) unitary operator $U:\ (\mathbb{C}^p\otimes\mathcal{M})\oplus\mathcal{E}\to (\mathbb{C}^q\otimes\mathcal{M})\oplus\mathcal{E}_*$ such that $U|\mbox{clos}(\mathcal{D}_0)=\widetilde{U_0}$.
In the case when \eqref{dim} doesn't hold one can set $\widetilde{\mathcal{M}}:=\mathcal{M\oplus K}$, where $\mathcal{K}$
is an infinite dimensional Hilbert space, then \eqref{dim} holds for $\widetilde{\mathcal{M}}$ in the place of $\mathcal{M}$, and there exists a unitary operator $U:\ (\mathbb{C}^p\otimes\widetilde{\mathcal{M}})\oplus\mathcal{E}\to (\mathbb{C}^q\otimes\widetilde{\mathcal{M}})\oplus\mathcal{E}_*$ such that $U|\mbox{clos}(\mathcal{D}_0)=\widetilde{U_0}$. Thus, without loss of generality we may consider that \eqref{dim} holds.

Let $U$ have a block partitioning
$$U=\left[\begin{array}{cc}
A & B \\
C & D
\end{array}\right]:\ (\mathbb{C}^p\otimes\mathcal{M})\oplus\mathcal{E}\to (\mathbb{C}^q\otimes\mathcal{M})\oplus\mathcal{E}_*.
$$
Then, in particular,
\begin{equation}\label{block}
\left[\begin{array}{cc}
A & B \\
C & D
\end{array}\right]\left[\begin{array}{c}
(P(w)\otimes I_\mathcal{M})H^R(w)\\
I_\mathcal{E}\end{array}\right] =\left[\begin{array}{c}
 H^R(w)\\
 \mathcal{F}(w)
\end{array}\right],\quad w\in\mathcal{D}_P.
\end{equation}
Since for $w\in\mathcal{D}_P$ one has $\| P(w)\| <1$, and since $\| A\|\leq 1$, we can solve the first block row equation of \eqref{block} for $H^R(w)$:
$$H^R(w)=(I_{\mathbb{C}^q\otimes\mathcal{H}}-A(P(w)\otimes I_\mathcal{M}))^{-1}B,\quad w\in\mathcal{D}_P.$$
Then from the second block row of \eqref{block} we get
$$\mathcal{F}(w)=D+C(P(w)\otimes I_\mathcal{M})(I_{\mathbb{C}^q\otimes\mathcal{M}}-A(P(w)\otimes I_\mathcal{M}))^{-1}B,\quad w\in\mathcal{D}_P,$$
i.e., \eqref{bb-tf} with $\mathcal{X}=\mathcal{M}$.

We are interested here in the case of the Agler--Schur class for the domain $\mathcal{D}_P$ where the domain $\mathcal{D}_P$ is $\mathbb{D}^{n_1\times n_1}\times\cdots\times\mathbb{D}^{n_N\times n_N}$, and the polynomial which defines this domain is $$P(W)=\mbox{diag}(W_1,\ldots, W_N),\quad W\in\mathbb{D}^{n_1\times n_1}\times\cdots\times\mathbb{D}^{n_N\times n_N}.$$
Here $W$ may be viewed as an ($n_1^2\cdots n_N^2$)-tuple of scalar variables $(W_k)_{ij},\ k=1,\ldots,N,\ i,j=1,\ldots,n_k$. We will write in this case $\mathcal{AS}_{n_1,\ldots,n_N}(\mathcal{E,E}_*)$ instead of $\mathcal{AS}_{\mathcal{D}_P}(\mathcal{E,E}_*)$, and if $\mathcal{E}=\mathcal{E}_*$ we will write $\mathcal{AS}_{n_1,\ldots,n_N}(\mathcal{E})$. The class $\mathcal{C}_{\mathcal{D}_P}$ is identified for $\mathcal{D}_P =\mathbb{D}^{n_1\times n_1}\times\cdots\times\mathbb{D}^{n_N\times n_N}$ with the class $\mathcal{C}^{(n_1,\ldots,n_N)}$ of $N$-tuples of matrices $\mathbf{T}=(T_1,\ldots,T_N)\in\mathcal{B}_\mathbf{T}^{n_1\times n_1}\times \cdots\times\mathcal{B}_\mathbf{T}^{n_N\times n_N}$ over a common commutative operator algebra $\mathcal{B}_\mathbf{T}\subset L(\mathcal{H}_\mathbf{T})$, with a Hilbert space $\mathcal{H}_\mathbf{T}$, such that $\| T_k\| <1,\ k=1,\ldots,N$.

Denote by $\mathcal{A}^{(n_1,\ldots,n_N)}$ the class of $N$-tuples of matrices $\mathbf{R}=(R_1,\ldots,R_N)\in\mathcal{B}_\mathbf{R}^{n_1\times n_1}\times \cdots\times\mathcal{B}_\mathbf{R}^{n_N\times n_N}$ over a common commutative operator algebra $\mathcal{B}_\mathbf{R}\subset L(\mathcal{H}_\mathbf{R})$, with a Hilbert space $\mathcal{H}_\mathbf{R}$, for which there exists a real constant $s_\mathbf{R}>0$ such that $$R_k+R_k^*\geq s_\mathbf{R}I_{\mathbb{C}_{n_k}\otimes\mathcal{H}_\mathbf{R}},\quad k=1,\ldots,N.$$
\begin{thm}\label{thm:spectr}
For any $\mathbf{R}\in\mathcal{A}^{(n_1,\ldots,n_N)}$,
$$\sigma_T(\mathbf{R})\subset\Pi^{n_1\times n_1}\times\cdots\times\Pi^{n_N\times n_N},$$
where $\sigma_T(\mathbf{R})$ denotes the Taylor joint spectrum of the collection of operators
$(R_k)_{ij},\ k=1,\ldots,N,\ i,j=1,\ldots,n_k$.
\end{thm}
\begin{proof}
It is shown in \cite{SZ} that the Taylor joint spectrum $\sigma_T(\mathbf{X})$ of an $n$-tuple of commuting bounded operators $\mathbf{X}=(X_1,\ldots,X_n)$ on a common Hilbert space $\mathcal{H}_\mathbf{X}$ is contained in the polynomially convex closure of $\sigma_\pi(\mathbf{X})$, the \emph{approximate point spectrum of $\mathbf{X}$}. The latter is defined as the set of points $\lambda=(\lambda_1,\ldots,\lambda_n)\in\mathbb{C}^n$ for which there exists a sequence of vectors $h_\nu\in\mathcal{H}_\mathbf{X}$ such that $\| h_\nu\| =1,\ \nu\in\mathbb{N}$, and $(X_j-\lambda_jI_{\mathcal{H}_\mathbf{X}})h_\nu\rightarrow 0$ as $\nu\rightarrow\infty$ for all $j=1,\ldots,n$. Thus it suffices to show that $$\sigma_\pi(\mathbf{R}):=\sigma_\pi\left(\{ (\mathbf{R}_k)_{ij}:\ k=1,\ldots,N,\ i,j=1,\ldots,n_k\}\right)\subset\Pi_s^{n_1\times n_1}\times\cdots\times\Pi_s^{n_N\times n_N}$$
whenever $\mathbf{R}\in\mathcal{A}^{(n_1,\ldots,n_N)}$ and $R_k+R_k^*\geq sI_{\mathbb{C}^{n_k}\otimes\mathcal{H}_\mathbf{R}}>0,\ k=1,\ldots,N$, where
$$\Pi_s^{n\times n}:=\{ M\in\mathbb{C}^{n\times n}: M+M^*\geq sI_n\},$$
since $\Pi_s^{n_1\times n_1}\times\cdots\times\Pi_s^{n_N\times n_N}$ is convex, and hence polynomially convex, and since $\Pi_s^{n_1\times n_1}\times\cdots\times\Pi_s^{n_N\times n_N}\subset\Pi^{n_1\times n_1}\times\cdots\times\Pi^{n_N\times n_N}$ for $s>0$. Suppose that $\Lambda=(\Lambda_1,\ldots,\Lambda_N)\in\sigma_\pi(\mathbf{R})$. Then there exists a sequence of vectors $h_\nu\in\mathcal{H}_\mathbf{R}$ such that $\| h_\nu\| =1,\ \nu\in\mathbb{N}$, and for $k=1,\ldots,N,\ i,j=1,\ldots,n_k$  one has
$$((R_k)_{ij}-(\Lambda_k)_{ij}I_{\mathcal{H}_\mathbf{R}})h_\nu\rightarrow 0\quad\mbox{as}\  \nu\rightarrow\infty.$$
Therefore, for every $k\in\{ 1,\ldots,N\}$ and $u_k=\mbox{col}(u_{k1},\ldots,u_{kn_k})\in\mathbb{C}^{n_k}$ one has
$$\sum\limits_{i=1}^{n_k}\sum_{j=1}^{n_k}\left(\left\langle ((R_k)_{ij}+(R_k)_{ji}^*)h_\nu,h_\nu\right\rangle -((\Lambda_k)_{ij}+\overline{(\Lambda_k)_{ji}})\left\langle h_\nu,h_\nu\right\rangle\right) u_{ki}\overline{u_{kj}}\rightarrow 0$$
as $\nu\rightarrow\infty$. Since $\left\langle h_\nu,h_\nu\right\rangle=1$, the subtrahend does not depend on $\nu$. Therefore,
\begin{eqnarray*}
s\left\langle u_k,u_k\right\rangle &=& s\lim_{\nu\rightarrow\infty}\left\langle u_k\otimes h_\nu,u_k\otimes h_\nu \right\rangle\\
&\leq & \lim_{\nu\rightarrow\infty}\left\langle \left( R_k+R_k^*\right)u_k\otimes h_\nu,u_k\otimes h_\nu\right\rangle \\ &=&\lim_{\nu\rightarrow\infty}\sum\limits_{i=1}^{n_k}\sum_{j=1}^{n_k}\left\langle ((R_k)_{ij}+(R_k)_{ji}^*)h_\nu,h_\nu\right\rangle u_{ki}\overline{u_{kj}} \\
&=& \sum\limits_{i=1}^{n_k}\sum_{j=1}^{n_k}((\Lambda_k)_{ij}+\overline{(\Lambda_k)_{ji}}) u_{ki}\overline{u_{kj}}
= \left\langle (\Lambda_k+\Lambda_k^*)u_k,u_k\right\rangle.
\end{eqnarray*}
Thus, $\Lambda_k+\Lambda_k^*\geq sI_{n_k},\ k=1,\ldots,N$, i.e., $\Lambda\in \Pi_s^{n_1\times n_1}\times\cdots\times\Pi_s^{n_N\times n_N}$, as desired.
\end{proof}
Theorem~\ref{thm:spectr} implies that for every holomorphic function $f$ on $\Pi^{n_1\times n_1}\times\cdots\times\Pi^{n_N\times n_N}$ and every $\mathbf{R}\in\mathcal{A}^{(n_1,\ldots,n_N)}$ the operator $f(\mathbf{R})$ is well defined by the Taylor functional calculus.

The \emph{Cayley transform} defined by
\begin{equation}\label{cal1}
R_k=(I_{\mathbb{C}^{n_k}\otimes\mathcal{H}_\mathbf{T}}+T_k)(I_{\mathbb{C}^{n_k}\otimes\mathcal{H}_\mathbf{T}}-T_k)^{-1},\quad k=1,\ldots,N,
\end{equation}
maps the class $\mathcal{C}^{(n_1,\ldots,n_N)}$ onto the class $\mathcal{A}^{(n_1,\ldots,n_N)}$, and its inverse is given by
\begin{equation}\label{cal2}
T_k=(R_k-I_{\mathbb{C}^{n_k}\otimes\mathcal{H}_\mathbf{R}})(R_k+I_{\mathbb{C}^{n_k}\otimes\mathcal{H}_\mathbf{R}})^{-1},\quad k=1,\ldots,N,
\end{equation}
where $\mathcal{H}_\mathbf{R}=\mathcal{H}_\mathbf{T}$. Let $f$ be an $L(\mathcal{U})$-valued function holomorphic on $\Pi^{n_1\times n_1}\times\cdots\times\Pi^{n_N\times n_N}$. Then its double Cayley transform $\mathcal{F}=\mathcal{C}(f)$ defined by \eqref{double-C-matr} is holomorphic on $\mathbb{D}^{n_1\times n_1}\times\cdots\times\mathbb{D}^{n_N\times n_N}$, and by the spectral mapping theorem and uniqueness of Taylor's functional calculus  (see \cite{P}) one has $$\mathcal{F}(\mathbf{T})=f(\mathbf{R}),$$ where $\mathbf{T}\in\mathcal{C}^{(n_1,\ldots,n_N)}$ and $\mathbf{R}\in\mathcal{A}^{(n_1,\ldots,n_N)}$ are related by \eqref{cal1} and \eqref{cal2}.
\begin{thm}\label{thm:b-calc}
Let $f$ be an $L(\mathcal{U})$-valued function holomorphic on $\Pi^{n_1\times n_1}\times\cdots\times\Pi^{n_N\times n_N}$. Then $f\in\mathcal{B}_{n_1,\ldots,n_N}(\mathcal{U})$ if and only if the following conditions are satisfied:
\begin{description}
  \item[(i)] $f(tZ)=tf(Z),\quad t>0,\
   Z\in\Pi^{n_1\times n_1}\times\cdots\times\Pi^{n_N\times n_N}$;
    \item[(ii)] $f(\mathbf{R})+f(\mathbf{R})^*\geq 0,\quad
    \mathbf{R}\in\mathcal{A}^{(n_1,\ldots,n_N)}$;
    \item[(iii)] $f(Z^*):=f(Z_1^*,\ldots,Z_N^*)=f(Z)^*,\quad  Z\in\Pi^{n_1\times n_1}\times\cdots\times\Pi^{n_N\times n_N}$.
\end{description}
\end{thm}
\begin{proof}
\textbf{Necessity.} Let $f\in\mathcal{B}_{n_1,\ldots,n_N}(\mathcal{U})$. Then
(i) and (iii) easily follow from the representation \eqref{sc}
of $f$. Condition (ii) on $f$ is equivalent to condition
\eqref{vN} on $\mathcal{F}$ which is defined by \eqref{double-C-matr}, i.e., to
$\mathcal{F}\in\mathcal{AS}_{n_1,\ldots,n_N}(\mathcal{U})$. Let us show the latter. Since by Theorem~\ref{thm:b-decomp}
$f$ satisfies \eqref{b-decomp}, and hence \eqref{b-decomp-1}, one can set
$$Z_k=(I_{n_k}+W_k)(I_{n_k}-W_k)^{-1},\
 \Lambda_k=(I_{n_k}+\Xi_k)(I_{n_k}-\Xi_k)^{-1},\ k=1,\ldots,N,$$ in
\eqref{b-decomp-1} and get
\begin{eqnarray*}
\lefteqn{(I_\mathcal{U}+\mathcal{F}(W))(I_\mathcal{U}-\mathcal{F}(W))^{-1}+(I_\mathcal{U}-\mathcal{F}(\Xi)^*)^{-1}(I_\mathcal{U}+\mathcal{F}(\Xi)^*) }\\
& =& \sum\limits_{k=1}^N\theta^\circ_k(\Xi)^*\left\{\left((I_{n_k}+W_k)(I_{n_k}-W_k)^{-1}+(I_{n_k}-\Xi_k^*)^{-1}(I_{n_k}+\Xi_k^*)\right)\otimes I_{\mathcal{M}_k}\right\}\\
& \times &\theta^\circ_k(W),\quad W,\Xi\in\mathbb{D}^{n_1\times n_1}\times\cdots\times\mathbb{D}^{n_N\times n_N},
\end{eqnarray*}
where for $k=1,\ldots,N$,
\begin{equation} \label{theta-circ}
    \theta^\circ_k(W)=\varphi_k\left((I_{n_1}+W_1)(I_{n_1}-W_1)^{-1},\ldots,(I_{n_N}+W_N)(I_{n_N}-W_N)^{-1}\right).
\end{equation}
We can rewrite this in the form
\begin{equation}\label{bb-id}
 I_\mathcal{U}-\mathcal{F}(\Xi)^*\mathcal{F}(W)=  \sum_{k=1}^N\theta_k(\Xi)^*\left((I_{n_k}-\Xi_k^*W_k)\otimes I_{\mathcal{M}_k}\right)\theta_k(W),
\end{equation}
where for $k=1,\ldots,N$,
\begin{equation}\label{theta}
    \theta_k(W)=\left((I_{n_k}-W_k)^{-1}\otimes I_{\mathcal{M}_k}\right)\theta^\circ_k(W)(I_\mathcal{U}-\mathcal{F}(W))\in L(\mathcal{U},\mathbb{C}^{n_k}\otimes\mathcal{M}_k).
\end{equation}
The identity \eqref{bb-id} coincides with \eqref{bb-id-r} for our case, with
$$H^R(W)=\mbox{col}(\theta_1(W),\ldots,\theta_N(W))\in L\left(\mathcal{U},\bigoplus\limits_{k=1}^N\left(\mathbb{C}^{n_k}\otimes\mathcal{M}_k\right)\right),$$
$$P(W)=\mbox{diag}(W_1,\ldots,W_N).$$
Note, that without loss of generality we may consider all of $\mathcal{M}_k$'s equal to some space $\mathcal{M}$, say, $\mathcal{M}=\bigoplus_{k=1}^N\mathcal{M}_k$. Then $H^R(W)\in L\left(\mathcal{U},\mathbb{C}^{n_1+\cdots +n_N}\otimes\mathcal{M}\right)$. By Theorem~\ref{thm:BB}, this means that $\mathcal{F}\in\mathcal{AS}_{n_1,\ldots,n_N}(\mathcal{U})$.

\textbf{Sufficiency.} Let $f$ satisfy conditions (i)--(iii). Since
(ii) is equivalent to $\mathcal{F}\in\mathcal{AS}_{n_1,\ldots,n_N}(\mathcal{U})$, where $\mathcal{F}$
is defined by \eqref{double-C-matr}, the identity \eqref{bb-id} holds
with some $L\left(\mathcal{U},\mathbb{C}^{n_k}\otimes\mathcal{M}\right)$-valued functions $\theta_k$ holomorphic  on
$\mathbb{D}^{n_1\times n_1}\times\cdots\times\mathbb{D}^{n_N\times n_N},\ k=1,\ldots,N$, with an auxiliary Hilbert space $\mathcal{M}$ (spaces $\mathcal{M}_k$ can be chosen equal in \eqref{bb-id}). Set
$$W_k=(Z_k-I_{n_k})(Z_k+I_{n_k})^{-1},\
\Xi_k=( \Lambda_k-I_{n_k})( \Lambda_k+I_{n_k})^{-1},\ k=1,\ldots,N,$$ in
\eqref{bb-id}, and by virtue of \eqref{double-C-matr} get \eqref{b-decomp-1} with
\begin{eqnarray}
\varphi_k(Z) &=& ((I_{n_k}+Z_k)^{-1})\otimes I_\mathcal{M})\nonumber \\
&\times & \theta_k\left((Z_1-I_{n_1})(Z_1+I_{n_1})^{-1},\ldots,(Z_N-I_{n_N})(Z_N+I_{n_N})^{-1}\right)\nonumber\\
&\times & (I_\mathcal{U}+f(Z))\in L(\mathcal{U},\mathbb{C}^{n_k}\otimes\mathcal{M}),\quad k=1,\ldots,N \label{fi}
\end{eqnarray}
(in fact, passing from \eqref{b-decomp-1} to \eqref{bb-id} is invertible, and \eqref{fi} is obtained from \eqref{theta-circ} and \eqref{theta}, and vice versa).
 The property (iii) implies $f(X)=f(X)^*$ for
every $N$-tuple $X=(X_1,\ldots,X_N)\in\Pi^{n_1\times n_1}\times\cdots\times\Pi^{n_N\times n_N}$ of positive definite matrices (we will denote this set by $\mathcal{P}^{(n_1,\ldots,n_N)}$), and for any such $X$ and
$t>0$ by \eqref{b-decomp-1} one has:
\begin{eqnarray*}
f(X)+f(tX) &=& (1+t)\sum_{k=1}^N\varphi_k(tX)^*(X_k\otimes I_\mathcal{M})\varphi_k(X),\\
f(tX)+f(X) &=& (1+t)\sum_{k=1}^N\varphi_k(X)^*(X_k\otimes I_\mathcal{M})\varphi_k(tX),\\
\frac{1+t}{2}[f(X)+f(X)] &=& \frac{1+t}{2}\sum_{k=1}^N\varphi_k(X)^*(2X_k\otimes I_\mathcal{M})\varphi_k(X),\\
\frac{1+t}{2t}[f(tX)+f(tX)] &=& \frac{1+t}{2t}
\sum_{k=1}^N2\varphi_k(tX)^*(2tX_k\otimes I_\mathcal{M})\varphi_k(tX).
\end{eqnarray*}
By (i), the left-hand sides of these equalities coincide and equal
 $(1+t)f(X)$, hence
\begin{eqnarray*}
f(X) &=& \sum_{k=1}^N\varphi_k(tX)^*(X_k\otimes I_\mathcal{M})\varphi_k(X)=
\sum_{k=1}^N\varphi_k(X)^*(X_k\otimes I_\mathcal{M})\varphi_k(tX)\\
     &=& \sum_{k=1}^N\varphi_k(X)^*(X_k\otimes I_\mathcal{M})\varphi_k(X)=\sum_{k=1}^N\varphi_k(tX)^*(X_k\otimes I_\mathcal{M})\varphi_k(tX).
\end{eqnarray*}
It follows from the latter equalities that
\begin{eqnarray*}
0 &\leq &
\sum_{k=1}^N[\varphi_k(tX)-\varphi_k(X)]^*(X_k\otimes I_\mathcal{M})[\varphi_k(tX)-\varphi_k(X)]\\
  &=& \sum_{k=1}^N\varphi_k(tX)^*(X_k\otimes I_\mathcal{M})\varphi_k(tX)-
  \sum_{k=1}^N\varphi_k(tX)^*(X_k\otimes I_\mathcal{M})\varphi_k(X)\\
  &-& \sum_{k=1}^N\varphi_k(X)^*(X_k\otimes I_\mathcal{M})\varphi_k(tX)+
  \sum_{k=1}^N\varphi_k(X)^*(X_k\otimes I_\mathcal{M})\varphi_k(X)=0.
\end{eqnarray*}
Thus $\varphi_k(tX)-\varphi_k(X)=0$ for every $X\in\mathcal{P}^{(n_1,\ldots,n_N)}$, $t>0$
and $k=1,\ldots,N$. For fixed
$k\in\{ 1,\ldots,N\}$ and $t>0$ the function
$h_{k,t}(Z):=\varphi_k(tZ)-\varphi_k(Z)$ is holomorphic on $\Pi^{n_1\times n_1}\times\cdots\times\Pi^{n_N\times n_N}$
and takes values in $L(\mathcal{U},\mathbb{C}^{n_k}\otimes\mathcal{M})$. Then for any fixed
$k\in\{ 1,\ldots,N\},\ t>0,\ u\in\mathcal{U}$ and
$m\in\mathbb{C}^{n_k}\otimes\mathcal{M}$ the scalar function $h_{k,t,u,m}(Z):=\langle
h_{k,t}(Z)u,m\rangle_{\mathbb{C}^{n_k}\otimes\mathcal{M}}$ is holomorphic on $\Pi^{n_1\times n_1}\times\cdots\times\Pi^{n_N\times n_N}$
and vanishes on $\mathcal{P}^{(n_1,\ldots,n_N)}$. The latter set is the uniqueness subset in $\Pi^{n_1\times n_1}\times\cdots\times\Pi^{n_N\times n_N}$, thus by the uniqueness theorem
for holomorphic functions of several variables (see, e.g.,
\cite{Sh2}), $h_{k,t,u,m}(Z)\equiv 0$, hence $h_{k,t}(Z)\equiv 0$,
which means:
$$\varphi_k(tZ)=\varphi_k(Z),\quad t>0,\ Z\in\Pi^{n_1\times n_1}\times\cdots\times\Pi^{n_N\times n_N}.$$
It follows from the latter equality that for every
$Z,\Lambda\in\Pi^{n_1\times n_1}\times\cdots\times\Pi^{n_N\times n_N}$ and $t>0$ one has
\begin{eqnarray*}
f(Z)+tf(\Lambda)^* &=& f(Z)+f(t\Lambda)^*=
\sum_{k=1}^N\varphi_k(t\Lambda)^*\left((Z_k+t\Lambda_k^*)\otimes I_\mathcal{M}\right)\varphi_k(Z)\\
 &=&
 \sum_{k=1}^N\varphi_k(\Lambda)^*\left((Z_k+t\Lambda_k^*)\otimes I_\mathcal{M}\right)\varphi_k(Z)\\
 &=&
\sum_{k=1}^N\varphi_k(\Lambda)^*\left(Z_k\otimes I_\mathcal{M}\right)\varphi_k(Z)+t\varphi_k(\Lambda)^*\left(\Lambda_k^*\otimes I_\mathcal{M}\right)\varphi_k(Z),
\end{eqnarray*}
and the comparison of the coefficients of the two linear functions
in $t$, at the beginning and at the end of this chain of
equalities, gives:
$$f(Z)=\sum_{k=1}^N\varphi_k(\Lambda)^*\left(Z_k\otimes I_\mathcal{M}\right)\varphi_k(Z),\quad Z,\Lambda\in\Pi^{n_1\times n_1}\times\cdots\times\Pi^{n_N\times n_N},$$ i.e., \eqref{b-decomp} with $\mathcal{M}_k=\mathcal{M},\ k=1,\ldots,N$. By
Theorem~\ref{thm:b-decomp}, $f\in\mathcal{B}_{n_1,\ldots,n_N}(\mathcal{U})$. The proof
is complete.
\end{proof}
\begin{cor}\label{cor:b-calc}
Let $f$ be an $L(\mathcal{U})$-valued function holomorphic on $\Omega_{n_1,\ldots,n_N}$. Then $f\in\mathcal{B}_{n_1,\ldots,n_N}(\mathcal{U})$ if and only if the following conditions are satisfied:
\begin{description}
    \item[(i)] $f(\lambda Z)=\lambda f(Z),\quad \lambda\in\mathbb{C}\setminus\{ 0\},\ Z\in \Omega_{n_1,\ldots,n_N}$;
    \item[(ii)] $f(\mathbf{R})+f(\mathbf{R})^*\geq 0,\quad \mathbf{R}\in\mathcal{A}^{(n_1,\ldots,n_N)}$;
    \item[(iii)] $f(Z^*)=f(Z)^*,\quad Z\in \Omega_{n_1,\ldots,n_N}$.
\end{description}
\end{cor}
\begin{proof}
If $f\in\mathcal{B}_{n_1,\ldots,n_N}(\mathcal{U})$ then  (i) and (iii) follow from the representation \eqref{sc}--\eqref{lp} of $f$, and (ii) follows from Theorem~\ref{thm:b-calc}.

Conversely, statements (i)--(iii) of the corrollary imply statements (i)--(iii) of Theorem~\ref{thm:b-calc}, which in turn imply that $f\in\mathcal{B}_{n_1,\ldots,n_N}(\mathcal{U})$.
\end{proof}
\begin{rem}\label{rem:b-calc}
By Corollary~\ref{cor:b-calc}, its conditions (i)--(iii) on
holomorphic $L\mathcal{(U)}$-valued functions on $\Omega_{n_1,\ldots,n_N}$ give
an equivalent definition of the class
$\mathcal{B}_{n_1,\ldots,n_N}(\mathcal{U})$, which seems to be more natural than
the original definition given above in
``existence" terms.
\end{rem}

\section{The image of the class $\mathcal{B}_{n_1,\ldots,n_N}(\mathcal{U})$ under the double Cayley transform}\label{sec:image}
 It was shown in the proof of Theorem~\ref{thm:b-calc} that if $f\in\mathcal{B}_{n_1,\ldots,n_N}(\mathcal{U})$ then the double Cayley transform of $f$, $\mathcal{F}=\mathcal{C}(f)$, defined by \eqref{double-C-matr}, belongs to the Agler--Schur class $\mathcal{AS}_{n_1,\ldots,n_N}(\mathcal{U})$. In fact, we are able to proof a stronger statement.
\begin{thm}\label{thm:image}
A holomorphic $L\mathcal{(U)}$-valued function $\mathcal{F}$ on
$\mathbb{D}^{n_1\times n_1}\times\cdots\times\mathbb{D}^{n_N\times n_N}$ can be represented as $\mathcal{F}=\mathcal{C}(f)$
for some $f\in\mathcal{B}_{n_1,\ldots,n_N}(\mathcal{U})$ if and only if the
following conditions are fulfilled:
\begin{description}
    \item[(i)] There exist  a Hilbert space $\mathcal{X}$ and an operator
    \begin{equation}\label{cal-u}
    U=\left[
\begin{array}{cc}
    A & B\\
    C & D
\end{array}\right]\in L((\mathbb{C}^{n_1+\cdots +n_N}\otimes\mathcal{X})\oplus\mathcal{U})
\end{equation}
such that for $W=(W_1,\ldots,W_N)\in\mathbb{D}^{n_1\times n_1}\times\cdots\times\mathbb{D}^{n_N\times n_N}$ one has
\begin{equation}\label{cal-tf}
\mathcal{F}(W)=D+C(P(W)\otimes I_{\mathcal{X}})(I_{\mathbb{C}^{n_1+\cdots +n_N}\otimes\mathcal{X}}-A(P(W)\otimes I_{\mathcal{X}}))^{-1}B,
\end{equation}
where $P(W)=\rm{diag}(W_1,\ldots,W_N)$ and $U=U^*=U^{-1}$.
    \item[(ii)] $1\notin\sigma(\mathcal{F}(0))$.
\end{description}
\end{thm}
\begin{proof}
\textbf{Necessity.} Let $f\in\mathcal{B}_{n_1,\ldots,n_N}(\mathcal{U})$. Then
\eqref{b-decomp-1} and \eqref{b-decomp-2} hold. As we have shown in
Theorem~\ref{thm:b-calc}, the identity \eqref{b-decomp-1} implies the
identity \eqref{bb-id} for
 $\mathcal{F}=\mathcal{C}(f)$, with holomorphic $L(\mathcal{U},\mathbb{C}^{n_k}\otimes\mathcal{M}_k)$-valued functions
$\theta_k,\ k=1,\ldots,N$, on
$\mathbb{D}^{n_1\times n_1}\times\cdots\times\mathbb{D}^{n_N\times n_N}$ defined by \eqref{theta-circ} and \eqref{theta}.
Analogously, the identity \eqref{b-decomp-2} implies
\begin{eqnarray}
\mathcal{F}(W)-\mathcal{F}(\Xi)^*
&= & \sum_{k=1}^N\theta_k(\Xi)^*\left((W_k-\Xi_k^*)\otimes I_{\mathcal{M}_k}\right)
\theta_k(W),\label{bb-id-2}\\
& & W,\Xi\in\mathbb{D}^{n_1\times n_1}\times\cdots\times\mathbb{D}^{n_N\times n_N}.\nonumber
\end{eqnarray}
Let us rewrite \eqref{bb-id} and \eqref{bb-id-2} in a somewhat different way. Since by Theorem~\ref{thm:b-calc} $f\in\mathcal{B}_{n_1,\ldots,n_N}(\mathcal{U})$ satisfies  $f(Z^*)=f(Z)^*,\ Z\in\Pi^{n_1\times n_1}\times\cdots\times\Pi^{n_N\times n_N}$, one has also $$\mathcal{F}(W^*)=\mathcal{F}(W)^*,\quad W\in\mathbb{D}^{n_1\times n_1}\times\cdots\times\mathbb{D}^{n_N\times n_N}.$$
Therefore, \eqref{bb-id} and \eqref{bb-id-2} are equivalent to the following two identities, respectively:
\begin{eqnarray}
I_\mathcal{U}-\mathcal{F}(W)\mathcal{F}(\Xi)^*
&= & \sum_{k=1}^N\widetilde{\theta_k}(W)\left((I_{n_k}-W_k\Xi_k^*)\otimes I_{\mathcal{M}_k}\right)
\widetilde{\theta_k}(\Xi)^*,\label{bb-id-1'}\\
\mathcal{F}(W)-\mathcal{F}(\Xi)^*
&= & \sum_{k=1}^N\widetilde{\theta_k}(W)\left((W_k-\Xi_k^*)\otimes I_{\mathcal{M}_k}\right)
\widetilde{\theta_k}(\Xi)^*,\label{bb-id-2'}\\
& & W,\Xi\in\mathbb{D}^{n_1\times n_1}\times\cdots\times\mathbb{D}^{n_N\times n_N},\nonumber
\end{eqnarray}
where $\widetilde{\theta_k}(W)={\theta_k(W^*)}^*$ are
holomorphic $L(\mathbb{C}^{n_k}\otimes\mathcal{M}_k,\mathcal{U})$-valued functions on
$\mathbb{D}^{n_1\times n_1}\times\cdots\times\mathbb{D}^{n_N\times n_N}$. We will show that the identities \eqref{bb-id-1'}
and \eqref{bb-id-2'} allow us to construct a Hilbert space $\mathcal{X}$ and an operator $U$ satisfying condition (i) of this theorem. To this end, we will apply the construction from Section~\ref{sec:b-calc} (next to Theorem~\ref{thm:BB}) to $\mathcal{F}=\mathcal{C}(f)$.  In this
case $\mathcal{E}=\mathcal{E}_*=\mathcal{U}$. Without loss of generality we may consider all of $\mathcal{M}_k$'s equal. Say, set $\mathcal{M}:=\bigoplus_{k=1}^N\mathcal{M}_k$ and regard $$H_k^R=\theta_k\in L(\mathcal{U},\mathbb{C}^{n_k}\otimes\mathcal{M}),\quad H_k^L=\widetilde{\theta_k}\in L(\mathbb{C}^{n_k}\otimes\mathcal{M},\mathcal{U}),\quad k=1,\ldots,N.$$ Then \eqref{bb-id}, \eqref{bb-id-2}, \eqref{bb-id-1'} and \eqref{bb-id-2'} imply \eqref{bb-id-lr}, and $$H^L(W^*)=H^R(W)^*,$$ where $H^L(W)=\rm{col}(H^L_1(W),\ldots,H^L_N(W)),\ H^R=\rm{col}(H^R_1(W),\ldots,H^R_N(W)),\ W\in\mathbb{D}^{n_1\times n_1}\times\cdots\times\mathbb{D}^{n_N\times n_N}$. Thus,
 $\mathcal{D}_0=\mathcal{R}_0$, and the operator $U_0$ acts on the generating vectors of $\mathcal{D}_0$ as follows:
 \begin{eqnarray*}
 \left[\begin{array}{c}
 (P(W)\otimes I_\mathcal{M})H^R(W)\\
 I_\mathcal{U}
 \end{array}\right]u &\longmapsto & \left[\begin{array}{c}
 H^R(W)\\
 \mathcal{F}(W)
 \end{array}\right]u,\\
 \left[\begin{array}{c}
 H^R(W)\\
 \mathcal{F}(W)
 \end{array}\right]u &\longmapsto &  \left[\begin{array}{c}
 (P(W)\otimes I_\mathcal{M})H^R(W)\\
 I_\mathcal{U}
 \end{array}\right]u,\\
 & & W\in\mathbb{D}^{n_1\times n_1}\times\cdots\times\mathbb{D}^{n_N\times n_N}.
  \end{eqnarray*}
 We used here the relations $\mathcal{F}(W^*)=\mathcal{F}(W)^*,\ H^L(W^*)=H^R(W)^*,\ P(W^*)=\rm{diag}(W_1^*,\ldots,W_N^*)=P(W)^*$. Thus  $U_0=U_0^{-1}$. Therefore,
$\widetilde{U_0}=\widetilde{U_0}^{-1}$. Since $p=q=n_1+\cdots +n_N,\ \mathcal{E}=\mathcal{E}_*=\mathcal{U}$, \eqref{dim} holds. Then the operator $$U=\widetilde{U_0}\oplus I_{((\mathbb{C}^{n_1+\cdots +n_N}\otimes\mathcal{M})\oplus\mathcal{U})\ominus\rm{clos}(\mathcal{D}_0)}\in L((\mathbb{C}^{n_1+\cdots +n_N}\otimes\mathcal{M})\oplus\mathcal{U})$$
satisfies $U=U^{-1}$. Since we have also $U^*=U^{-1}$, (i) is satisfied with $\mathcal{X}=\mathcal{M}$.

Statement (ii) follows in the same way as in \cite[Theorem 4.2]{K}, with $E=(I_{n_1},\ldots,I_{n_N})$ in the place of $e=(1,\ldots,1)$.

\textbf{Sufficiency.} Let the conditions (i) and (ii) on
$\mathcal{F}$ be satisfied. Then in the same way as in \cite[Theorem 4.2]{K} one can see that $1\notin\sigma(\mathcal{F}(W))$ for all $W\in\mathbb{D}^{n_1\times n_1}\times\cdots\times\mathbb{D}^{n_N\times n_N}$. Thus,  the function
$$F(W):=(I_\mathcal{U}+\mathcal{F}(W))(I_\mathcal{U}-\mathcal{F}(W))^{-1}$$
is correctly defined and holomorphic on $\mathbb{D}^{n_1\times n_1}\times\cdots\times\mathbb{D}^{n_N\times n_N}$. It is easy
to see that
\begin{equation}\label{+}
    F(W)+F(\Xi)^*=
    2(I_\mathcal{U}-\mathcal{F}(\Xi)^*)^{-1}
    (I_\mathcal{U}-\mathcal{F}(\Xi)^*\mathcal{F}(W))
    (I_\mathcal{U}-\mathcal{F}(W))^{-1},
\end{equation}
\begin{equation}\label{-}
    F(W)-F(\Xi)^*=
    2(I_\mathcal{U}-\mathcal{F}(\Xi)^*)^{-1}
    (\mathcal{F}(W)-\mathcal{F}(\Xi)^*)
    (I_\mathcal{U}-\mathcal{F}(W))^{-1}.
\end{equation}
As shown in \cite{BB}, it follows from \eqref{cal-tf} that
\begin{eqnarray*}I_\mathcal{U}-\mathcal{F}(\Xi)^*\mathcal{F}(W) &=&
 B^*(I_{\mathbb{C}^{n_1+\cdots +n_N}\otimes\mathcal{X}}-(P(\Xi)^*\otimes I_\mathcal{X})A^*)^{-1}\\
 &\times & ((I_{\mathbb{C}^{n_1+\cdots +n_N}}-P(\Xi)^*P(W))\otimes I_\mathcal{X})\\
 &\times & (I_{\mathbb{C}^{n_1+\cdots +n_N}\otimes\mathcal{X}}-A(P(W)\otimes I_\mathcal{X}))^{-1}B.
 \end{eqnarray*}
 Since $U=U^*$, we get
 \begin{eqnarray*}
 I_\mathcal{U}-\mathcal{F}(\Xi)^*\mathcal{F}(W) &=&
 \sum\limits_{k=1}^NB^*(I_{\mathbb{C}^{n_1+\cdots +n_N}\otimes\mathcal{X}}-(P(\Xi)^*\otimes I_\mathcal{X})A)^{-1}\\
 &\times & (P_{\mathbb{C}^{n_k}}\otimes I_\mathcal{X})((I_{n_k}-\Xi_k^*W_k)\otimes I_\mathcal{X})(P_{\mathbb{C}^{n_k}}\otimes I_\mathcal{X})\\
 &\times & (I_{\mathbb{C}^{n_1+\cdots +n_N}\otimes\mathcal{X}}-A(P(W)\otimes I_\mathcal{X}))^{-1}B.
 \end{eqnarray*}
 Analogously,
 \begin{eqnarray*}
\mathcal{F}(W)- \mathcal{F}(\Xi)^* &=&
 \sum\limits_{k=1}^NB^*(I_{\mathbb{C}^{n_1+\cdots +n_N}\otimes\mathcal{X}}-(P(\Xi)^*\otimes I_\mathcal{X})A)^{-1}\\
 &\times & (P_{\mathbb{C}^{n_k}}\otimes I_\mathcal{X})((W_k-\Xi_k^*)\otimes I_\mathcal{X})(P_{\mathbb{C}^{n_k}}\otimes I_\mathcal{X})\\
 &\times & (I_{\mathbb{C}^{n_1+\cdots +n_N}\otimes\mathcal{X}}-A(P(W)\otimes I_\mathcal{X}))^{-1}B.
 \end{eqnarray*}
 Thus, from \eqref{+} and \eqref{-} we get
 \begin{eqnarray}
  F(W)+F(\Xi)^* &=& \sum\limits_{k=1}^N\xi_k(\Xi)^*((I_{n_k}-\Xi_k^*W_k)\otimes I_\mathcal{X})\xi_k(W),\label{++}\\
  F(W)-F(\Xi)^* &=& \sum\limits_{k=1}^N\xi_k(\Xi)^*((W_k-\Xi_k^*)\otimes I_\mathcal{X})\xi_k(W),\label{--}
 \end{eqnarray}
 with
 $$\xi_k(W)=\sqrt{2}(P_{\mathbb{C}^{n_k}\otimes I_\mathcal{X}})(I_{\mathbb{C}^{n_1+\cdots +n_N}\otimes\mathcal{X}}-A(P(W)\otimes I_\mathcal{X}))^{-1}B(I_\mathcal{U}-\mathcal{F}(W))^{-1},$$
 for all $W\in\mathbb{D}^{n_1\times n_1}\times\cdots\times\mathbb{D}^{n_N\times n_N}$ and $k=1,\ldots,N$.
  Since for $Z_k,\Lambda_k\in\Pi^{n_k\times n_k}$ we have
  \begin{eqnarray*}
  \lefteqn{I_{n_k}-(\Lambda_k^*+I_{n_k})^{-1}(\Lambda_k^*-I_{n_k})(Z_k-I_{n_k})(Z_k+I_{n_k})^{-1}}\\
  &=2(\Lambda_k^*+I_{n_k})^{-1}(Z_k+\Lambda_k^*)(Z_k+I_{n_k})^{-1},\\
  \lefteqn{(Z_k-I_{n_k})(Z_k+I_{n_k})^{-1}-(\Lambda_k^*+I_{n_k})^{-1}(\Lambda_k^*-I_{n_k})}\\
  &=2(\Lambda_k^*+I_{n_k})^{-1}(Z_k-\Lambda_k^*)(Z_k+I_{n_k})^{-1},
  \end{eqnarray*}
 by setting $W_k:=(Z_k-I_{n_k})(Z_k+I_{n_k})^{-1}$ and $\Xi_k:=(\Lambda_k-I_{n_k})(\Lambda_k+I_{n_k})^{-1}$ in \eqref{++} and \eqref{--} we get the identities \eqref{b-decomp-1} and \eqref{b-decomp-2} for
 $$f(Z)=F((Z_1-I_{n_1})(Z_1+I_{n_1})^{-1},\ldots,(Z_N-I_{n_N})(Z_N-I_{n_N})^{-1}),$$ with
 \begin{eqnarray*}\varphi_k(Z) &=& \sqrt{2}((Z_k+I_{n_k})^{-1}\otimes I_{\mathcal X})\\
 &\times & \xi_k((Z_1-I_{n_1})(Z_1+I_{n_1})^{-1},\ldots,(Z_N-I_{n_N})(Z_N-I_{n_N})^{-1})
 \end{eqnarray*}
 for $k=1,\ldots,N$.
  Thus, by Theorem~\ref{thm:b-decomp} we finally get $\mathcal{F}=\mathcal{C}(f)$
where $f\in\mathcal{B}_{n_1,\ldots,n_N}(\mathcal{U})$. The proof is complete.
\end{proof}

\section{The ``real" case}\label{sec:real-b}
In Section~\ref{sec:intro} we have mentioned the notions of an anti-unitary involution (AUI) $\iota=\iota_\mathcal{U}$ on a Hilbert space $\mathcal{U}$ (a counterpart of the operator $\iota_n$ of complex conjugation on $\mathbb{C}^n$), a $(\iota_\mathcal{U},\iota_\mathcal{Y})$-real operator $A\in L(\mathcal{U},\mathcal{Y})$ (a counterpart of matrix with real entries), and a  $(\iota_\mathcal{U},\iota_\mathcal{Y})$-real operator-valued function $f$ (a counterpart of function which takes real scalar or matrix values   at real points). Some basic properties of AUI were described in \cite[Proposition 6.1]{K}. We will need also the following property.
\begin{prop}\label{prop:iota}
Let $\iota_\mathcal{U}$ and $\iota_\mathcal{H}$ be AUIs on Hilbert spaces $\mathcal{U}$ and $\mathcal{H}$, respectively. Then the operator $\iota_{\mathcal{U\otimes H}}=\iota_\mathcal{U}\otimes\iota_\mathcal{H}$ on $\mathcal{U\otimes H}$ which is defined on elementary tensors $u\otimes h$ as
\begin{equation}\label{iota-tensor}
(\iota_\mathcal{U}\otimes\iota_\mathcal{H})(u\otimes h)=\iota_\mathcal{U}u\otimes\iota_\mathcal{H}h
\end{equation}
and then extended to all of $\mathcal{U\otimes H}$ by linearity and continuity, is defined correctly and is an AUI on $\mathcal{U\otimes H}$.
\end{prop}
\begin{proof}
First, let us observe that $\iota_{\mathcal{U\otimes H}}$ is correctly defined. To this end, note that for arbitrary $x'=\sum_{\alpha=1}^lu'_\alpha\otimes h'_\alpha$ and $x''=\sum_{\beta=1}^mu''_\beta\otimes h''_\beta$ from $\mathcal{U\otimes H}$ we have
\begin{eqnarray*}
\lefteqn{\left\langle \iota_{\mathcal{U\otimes H}}x',\iota_{\mathcal{U\otimes H}}x''\right\rangle_{\mathcal{U\otimes H}} = \sum_{\alpha=1}^l\sum_{\beta=1}^m\left\langle\iota_\mathcal{U}u'_\alpha\otimes\iota_\mathcal{H}h'_\alpha,\iota_\mathcal{U}u''_\beta\otimes\iota_\mathcal{H}h''_\beta \right\rangle_{\mathcal{U\otimes H}}}\\
&=& \sum_{\alpha=1}^l\sum_{\beta=1}^m\left\langle\iota_\mathcal{U}u'_\alpha,\iota_\mathcal{U}u''_\beta\right\rangle_\mathcal{U}\left\langle \iota_\mathcal{H}h'_\alpha,\iota_\mathcal{H}h''_\beta \right\rangle_{\mathcal H}
= \sum_{\alpha=1}^l\sum_{\beta=1}^m\left\langle u''_\beta,u'_\alpha\right\rangle_\mathcal{U}\left\langle h''_\beta,h'_\alpha \right\rangle_{\mathcal H}\\
&=& \sum_{\alpha=1}^l\sum_{\beta=1}^m\left\langle u''_\beta\otimes h''_\beta,u'_\alpha\otimes h'_\alpha \right\rangle_{\mathcal{U\otimes H}}= \left\langle x'',x'\right\rangle_{\mathcal{U\otimes H}},
\end{eqnarray*}
i.e., $\iota_{\mathcal{U\otimes H}}$ is an \emph{anti-isometry} on linear combinations of elementary tensors. Thus, it is uniquely extended to an operator on all of $\mathcal{U\otimes H}$, and the property \eqref{antiun} of the extended operator follows by continuity. Since for arbitrary $x'=\sum_{\alpha=1}^lu'_\alpha\otimes h'_\alpha$ and $x''=\sum_{\beta=1}^mu''_\beta\otimes h''_\beta$ from $\mathcal{U\otimes H}$ we have
\begin{eqnarray*}
\left\langle \iota^2_{\mathcal{U\otimes H}}x',x''\right\rangle_{\mathcal{U\otimes H}} &=& \sum_{\alpha=1}^l\sum_{\beta=1}^m\left\langle\iota^2_\mathcal{U}u'_\alpha\otimes\iota^2_\mathcal{H}h'_\alpha,u''_\beta\otimes h''_\beta \right\rangle_{\mathcal{U\otimes H}}\\
&=& \sum_{\alpha=1}^l\sum_{\beta=1}^m\left\langle u'_\alpha\otimes h'_\alpha,u''_\beta\otimes h''_\beta \right\rangle_{\mathcal{U\otimes H}}\\
&=& \left\langle x',x''\right\rangle_{\mathcal{U\otimes H}},
\end{eqnarray*}
by continuity the property \eqref{invol}  of $\iota_{\mathcal{U\otimes H}}$ follows as well. Thus, $\iota_{\mathcal{U\otimes H}}$ is an AUI on $\mathcal{U\otimes H}$.
\end{proof}

Let $\mathcal{U}$ be a Hilbert space, and let $\iota=\iota_\mathcal{U}$ be an AUI on $\mathcal{U}$. Denote by $\iota\mathbb{R}\mathcal{B}_{n_1,\ldots,n_N}(\mathcal{U})$ the subclass of $\mathcal{B}_{n_1,\ldots,n_N}(\mathcal{U})$ consisting of $\iota$-real functions. The following theorem gives several equivalent characterizations of the ``$\iota$-real valued Bessmertny\u{\i} class" $\iota\mathbb{R}\mathcal{B}_{n_1,\ldots,n_N}(\mathcal{U})$ which specify for this case the characterizations obtained above for the ``complex valued Bessmertny\u{\i} class" $\mathcal{B}_{n_1,\ldots,n_N}(\mathcal{U})$.
\begin{thm}\label{thm:real}
Let $f$ be a holomorphic $L(\mathcal{U})$-valued function on
$\Omega_{n_1,\ldots,n_N}$, and $\iota=\iota_\mathcal{U}$ be an AUI on a Hilbert
space $\mathcal{U}$. The following statements are equivalent:
\begin{description}
    \item[(i)] $f\in\iota\mathbb{R}\mathcal{B}_{n_1,\ldots,n_N}(\mathcal{U})$;
    \item[(ii)] there exist a representation  \eqref{sc} of
    $f$ and AUIs $\iota_{\mathcal{M}_k}$ on $\mathcal{M}_k,\ k=1,\ldots,N$, such that
    the operators $G_k$ in  \eqref{lp} are
    $(\iota_\mathcal{U}\oplus\iota_\mathcal{H},\iota_{n_k}\otimes\iota_{\mathcal{M}_k})$-real;
    \item[(iii)] there exist a representation  \eqref{b-decomp}
    of $f$ and AUIs $\iota_{\mathcal{M}_k}$ on $\mathcal{M}_k$ such that the holomorphic functions  $\varphi_k(Z)$ are $(\iota_\mathcal{U},\iota_{n_k}\otimes\iota_{\mathcal{M}_k})$-real, $k=1,\ldots,N$;
    \item[(iv)] there exist a Hilbert space
$\mathcal{X}$ and an operator $U$ as in \eqref{cal-u} such that $\mathcal{F}=\mathcal{C}(f)$
 satisfies \eqref{cal-tf} and $U=U^*=U^{-1}$; moreover, there exists an AUI
    $\iota_{\mathcal{X}}$ on $\mathcal{X}$
    such that the operator $U$ is
    $((\iota_{n_1+\cdots +n_N}\otimes\iota_\mathcal{X})\oplus\iota_\mathcal{U})$-real.
\end{description}
\end{thm}
\begin{proof}
(i)$\Rightarrow $(iii) Let (i) hold. By Theorem~\ref{thm:b-decomp} there exists a representation \eqref{b-decomp} of $f$ with holomorphic $L(\mathcal{U},\mathbb{C}^{n_k}\otimes\mathcal{M}_k)$-valued functions $\varphi_k$ on $\Pi^{n_1\times n_1}\times\cdots\times\Pi^{n_N\times n_N}$. Let $\iota_{\mathcal{M}_k}$ be an AUI on $\mathcal{M}_k$, and let $\iota_{n_k}$ be a standard AUI on $\mathbb{C}^{n_k}$, i.e., a complex conjugation. Set $\widetilde{\mathcal{M}_k}:=\mathcal{M}_k\oplus\mathcal{M}_k$ and $\iota_{\widetilde{\mathcal{M}_k}}:=\left[\begin{array}{cc}
0 & \iota_{\mathcal{M}_k}\\
\iota_{\mathcal{M}_k} & 0
\end{array}\right],\ k=1,\ldots,N.$ Clearly, $\iota_{\widetilde{\mathcal{M}_k}}$ is an AUI on $\widetilde{\mathcal{M}_k}$. Define the rearrangement isomorphisms $V_k:\,\mathbb{C}^{n_k}\otimes (\mathcal{M}_k\oplus\mathcal{M}_k)\longrightarrow (\mathbb{C}^{n_k}\otimes \mathcal{M}_k)\oplus (\mathbb{C}^{n_k}\otimes \mathcal{M}_k)$ by
$$\left[\begin{array}{c}
m_{11}\\
m_{21}\\
\vdots\\
m_{1n_k}\\
m_{2n_k}\end{array}\right]\longmapsto\left[\begin{array}{c}
m_{11}\\
\vdots\\
m_{1n_k}\\
m_{21}\\
\vdots\\
m_{2n_k}\end{array}\right].$$
Then \begin{equation}\label{iota-iso}
\iota_{n_k}\otimes\iota_{\widetilde{\mathcal{M}_k}}=V_k^{-1}\left[\begin{array}{cc}
0 & \iota_{n_k}\otimes\iota_{\mathcal{M}_k}\\
\iota_{n_k}\otimes\iota_{\mathcal{M}_k} & 0
\end{array}\right]V_k.
\end{equation}
Set
$$\widetilde{\varphi_k}(Z):=\frac{1}{\sqrt{2}}V_k^{-1}\left[\begin{array}{c}
\varphi_k(Z)\\
(\iota_{n_k}\otimes\iota_{\mathcal{M}_k})\varphi_k(\bar{Z})\iota_\mathcal{U}\end{array}\right],$$
where $\bar{Z}=(\overline{Z_1},\ldots,\overline{Z_N})$, and $(\overline{Z_k})_{ij}=\overline{(Z_k)_{ij}},\ k=1,\ldots,N,\ i,j=1,\ldots,n_k$. By properties of AUIs,  $\widetilde{\varphi_k}(Z)$ is holomorphic on $\Pi^{n_1\times n_1}\times\cdots\times\Pi^{n_N\times n_N}$. Moreover, $\widetilde{\varphi_k}$ is $(\iota_\mathcal{U},\iota_{n_k}\otimes\iota_{\widetilde{\mathcal{M}_k}})$-real. Indeed, due to \eqref{iota-iso} we have
\begin{eqnarray*}
\widetilde{\varphi_k}^\sharp(Z) &=& (\iota_{n_k}\otimes\iota_{\widetilde{\mathcal{M}_k}})\widetilde{\varphi_k}(\bar{Z})\iota_\mathcal{U}\\ &=& \frac{1}{\sqrt{2}}V_k^{-1}\left[\begin{array}{cc}
0 & \iota_{n_k}\otimes\iota_{\mathcal{M}_k}\\
\iota_{n_k}\otimes\iota_{\mathcal{M}_k} & 0
\end{array}\right]\cdot\left[\begin{array}{c}
\varphi_k(\bar{Z})\\
(\iota_{n_k}\otimes\iota_{\mathcal{M}_k})\varphi_k(Z)\iota_\mathcal{U}\end{array}\right]\iota_\mathcal{U}\\
&=& \frac{1}{\sqrt{2}}V_k^{-1}\left[\begin{array}{c}
(\iota_{n_k}\otimes\iota_{\mathcal{M}_k})^2\varphi_k(Z)\iota_\mathcal{U}^2\\
(\iota_{n_k}\otimes\iota_{\mathcal{M}_k})\varphi_k(\bar{Z})\iota_\mathcal{U}\end{array}\right]\\
&=&
\frac{1}{\sqrt{2}}V_k^{-1}\left[\begin{array}{c}
\varphi_k(Z)\\
(\iota_{n_k}\otimes\iota_{\mathcal{M}_k})\varphi_k(\bar{Z})\iota_\mathcal{U}\end{array}\right]\\
&=& \widetilde{\varphi_k}(Z).
\end{eqnarray*}
Let us show that
$$f(Z)=\sum_{k=1}^N\widetilde{\varphi_k}(\Lambda)^*(Z_k\otimes I_{\widetilde{\mathcal{M}_k}})\widetilde{\varphi_k}(Z),\quad Z,\Lambda\in\Pi^{n_1\times n_1}\times\cdots\times\Pi^{n_N\times n_N}.$$
To this end, let us show first  that for $k=1,\ldots,N$:
\begin{equation}\label{*-id}
\left((\iota_{n_k}\otimes\iota_{\widetilde{\mathcal{M}_k}})\widetilde{\varphi_k}(\Lambda)\iota_\mathcal{U}\right)^*=\iota_\mathcal{U}\widetilde{\varphi_k}(\Lambda)^*(\iota_{n_k}\otimes\iota_{\widetilde{\mathcal{M}_k}}),\quad \Lambda\in\Pi^{n_1\times n_1}\times\cdots\times\Pi^{n_N\times n_N}.
\end{equation}
Indeed, for any $m\in\mathbb{C}^{n_k}\otimes\widetilde{\mathcal{M}_k},\ u\in\mathcal{U},\ \Lambda\in\Pi^{n_1\times n_1}\times\cdots\times\Pi^{n_N\times n_N}$ one has
\begin{eqnarray*}
\left\langle \left((\iota_{n_k}\otimes\iota_{\widetilde{\mathcal{M}_k}})\widetilde{\varphi_k}(\Lambda)\iota_\mathcal{U}\right)^*m,u\right\rangle_\mathcal{U} &=& \left\langle m,(\iota_{n_k}\otimes\iota_{\widetilde{\mathcal{M}_k}})\widetilde{\varphi_k}(\Lambda)\iota_\mathcal{U}u\right\rangle_{\mathbb{C}^{n_k}\otimes\widetilde{\mathcal{M}_k}}\\
\left\langle \widetilde{\varphi_k}(\Lambda)\iota_\mathcal{U}u,(\iota_{n_k}\otimes\iota_{\widetilde{\mathcal{M}_k}})m\right\rangle_{\mathbb{C}^{n_k}\otimes\widetilde{\mathcal{M}_k}}
&=& \left\langle \iota_\mathcal{U}u,\widetilde{\varphi_k}(\Lambda)^*(\iota_{n_k}\otimes\iota_{\widetilde{\mathcal{M}_k}})m\right\rangle_{\mathcal U}\\
&=& \left\langle \iota_\mathcal{U}\widetilde{\varphi_k}(\Lambda)^*(\iota_{n_k}\otimes\iota_{\widetilde{\mathcal{M}_k}})m,u\right\rangle_\mathcal{U}.
\end{eqnarray*}
Now, for any $Z,\Lambda\in\Pi^{n_1\times n_1}\times\cdots\times\Pi^{n_N\times n_N}$:
\begin{eqnarray*}
\lefteqn{\sum\limits_{k=1}^N\widetilde{\varphi_k}(\Lambda)^*(Z_k\otimes I_{\widetilde{\mathcal{M}_k}})\widetilde{\varphi_k}(Z)= \frac{1}{2}\sum\limits_{k=1}^N\left[\begin{array}{c}
\varphi_k(\Lambda)\\
(\iota_{n_k}\otimes\iota_{\mathcal{M}_k})\varphi_k(\bar{\Lambda})\iota_\mathcal{U}\end{array}\right]^*}\\
&\times & V_k(Z_k\otimes I_{\widetilde{\mathcal{M}_k}})V_k^{-1}\left[\begin{array}{c}
\varphi_k(Z)\\
(\iota_{n_k}\otimes\iota_{\mathcal{M}_k})\varphi_k(\bar{Z})\iota_\mathcal{U}\end{array}\right]\\
&=& \frac{1}{2}\sum\limits_{k=1}^N\left[\begin{array}{c}
\varphi_k(\Lambda)\\
(\iota_{n_k}\otimes\iota_{\mathcal{M}_k})\varphi_k(\bar{\Lambda})\iota_\mathcal{U}\end{array}\right]^*\left[\begin{array}{cc}
Z_k\otimes I_{\mathcal{M}_k} & 0\\
0 & Z_k\otimes I_{\mathcal{M}_k}
\end{array}\right]\\
&\times & \left[\begin{array}{c}
\varphi_k(Z)\\
(\iota_{n_k}\otimes\iota_{\mathcal{M}_k})\varphi_k(\bar{Z})\iota_\mathcal{U}\end{array}\right]\\
&=& \frac{1}{2}\sum\limits_{k=1}^N\varphi_k(\Lambda)^*( Z_k\otimes I_{\mathcal{M}_k})\varphi_k(Z)\\
&+& \frac{1}{2} \sum\limits_{k=1}^N\iota_\mathcal{U}\varphi_k(\bar{\Lambda})^*(\iota_{n_k}\otimes\iota_{\mathcal{M}_k})( Z_k\otimes I_{\mathcal{M}_k})(\iota_{n_k}\otimes\iota_{\mathcal{M}_k})\varphi_k(\bar{Z})\iota_\mathcal{U}\\
&=& \frac{1}{2}\left(\sum\limits_{k=1}^N\varphi_k(\Lambda)^*( Z_k\otimes I_{\mathcal{M}_k})\varphi_k(Z)+\sum\limits_{k=1}^N\iota_\mathcal{U}\varphi_k(\bar{\Lambda})^*(\overline{Z_k}\otimes I_{\mathcal{M}_k})\varphi_k(\bar{Z})\iota_\mathcal{U}\right)\\
&=& \frac{1}{2}(f(Z)+\iota_\mathcal{U}f(\bar{Z})\iota_\mathcal{U})=f(Z),
\end{eqnarray*}
where we used \eqref{*-id}, unitarity of $V_k$, and identity $\iota_{n_k}\overline{Z_k}\iota_{n_k}=Z_k$. Thus, (iii) follows from (i).

(iii)$\Rightarrow$(ii) Let (iii) hold. As in the sufficiency part of the proof of Theorem~\ref{b-decomp} we set
\begin{eqnarray*}
\mathcal{N}:=\bigoplus_{k=1}^N(\mathbb{C}^{n_k}\otimes\mathcal{M}_k), & \varphi(Z):={\rm col}(\varphi_1(Z),\ldots,\varphi_N(Z))\in L(\mathcal{U,N}),  \\ P_k:=P_{\mathcal{M}_k}, &  E=(I_{n_1},\ldots,I_{n_N})\in\Pi^{n_1\times n_1}\times\cdots\times\Pi^{n_N\times n_N},
\end{eqnarray*}
\begin{eqnarray*}  \mathcal{H}:=\rm{clos\ span}_{\Lambda\in\Pi^{n_1\times n_1}\times\cdots\times\Pi^{n_N\times n_N}}\left\{ (\varphi(\Lambda)-\varphi(E))\mathcal{U}\right\}\subset\mathcal{N}, \\  G_k:=(I_{n_k}\otimes P_k)\kappa\left[\begin{array}{cc}
\varphi(E) & 0\\
0 & I_\mathcal{H}\end{array}\right]\in L(\mathcal{U\oplus H},\mathbb{C}^{n_k}\otimes\mathcal{M}_k),
\end{eqnarray*}
where $\kappa:\,\mathcal{X\oplus H}\rightarrow\mathcal{N}$ is defined by \eqref{kappa}. For $\psi(Z)=\left[\begin{array}{c}
I_\mathcal{U}\\
\varphi(Z)-\varphi(E)
\end{array}\right]$ one has $\varphi(E)\mathcal{U}=\mathcal{U}\oplus\{ 0\}$, therefore the linear span of vectors of the form $\psi(Z)u,\ Z\in\Pi^{n_1\times n_1}\times\cdots\times\Pi^{n_N\times n_N},\ u\in\mathcal{U}$, is dense in $\mathcal{U\oplus H}$. Set $\iota_\mathcal{N}:=\bigoplus_{k=1}^N(\iota_{n_k}\otimes\iota_{\mathcal{M}_k})$. By the assumption, we have for $k=1,\ldots,N$: $$(\iota_{n_k}\otimes\iota_{\mathcal{M}_k})\varphi_k(Z)=\varphi_k(\bar{Z})\iota_\mathcal{U},\quad Z\in\Pi^{n_1\times n_1}\times\cdots\times\Pi^{n_N\times n_N}.$$ Therefore,
$$\iota_\mathcal{N}(\varphi(Z)-\varphi(E))u=(\varphi(Z)-\varphi(E))\iota_\mathcal{U}u\in\mathcal{H},\quad u\in\mathcal{U}.$$
Thus $\iota_\mathcal{N}\mathcal{H}\subset\mathcal{H}$. Moreover, $\mathcal{H}=\iota_\mathcal{N}^2\mathcal{H}\subset\iota_\mathcal{N}\mathcal{H}$, therefore $\iota_\mathcal{N}\mathcal{H}=\mathcal{H}$. Set $\iota_\mathcal{H}:=\iota_\mathcal{N}|\mathcal{H}$.
Clearly, $\iota_\mathcal{H}$ is an AUI on $\mathcal{H}$, and
$$(\iota_\mathcal{U}\oplus\iota_\mathcal{H})\psi(Z)=\psi(\bar{Z})\iota_\mathcal{U},\quad Z\in\Pi^{n_1\times n_1}\times\cdots\times\Pi^{n_N\times n_N}.$$
Let us verify that $G_k$ is $(\iota_\mathcal{U}\oplus\iota_\mathcal{H},\iota_{n_k}\otimes\iota_{\mathcal{M}_k})$-real, $k=1,\ldots,N$.
\begin{eqnarray*}
\lefteqn{(\iota_{n_k}\otimes\iota_{\mathcal{M}_k})G_k\psi(Z)u}\\
& =& (\iota_{n_k}\otimes\iota_{\mathcal{M}_k})(I_{n_k}\otimes P_k)\kappa\left[\begin{array}{cc}
\varphi(E) & 0\\
0 & I_\mathcal{H}
\end{array}\right]\left[\begin{array}{c}
I_\mathcal{U}\\
\varphi(Z)-\varphi(E)
\end{array}\right]u\\
&=& (\iota_{n_k}\otimes\iota_{\mathcal{M}_k})(I_{n_k}\otimes P_k)\kappa\left[\begin{array}{c}
\varphi(E)\\
\varphi(Z)-\varphi(E)
\end{array}\right]u = (\iota_{n_k}\otimes\iota_{\mathcal{M}_k})\varphi_k(Z)u\\
&=& \varphi_k(\bar{Z})\iota_\mathcal{U}u = (I_{n_k}\otimes P_k)\kappa\left[\begin{array}{cc}
\varphi(E) & 0\\
0 & I_\mathcal{H}
\end{array}\right]\left[\begin{array}{c}
I_\mathcal{U}\\
\varphi(\bar{Z})-\varphi(E)
\end{array}\right]\iota_\mathcal{U}u\\
&=& G_k\psi(\bar{Z})\iota_\mathcal{U}u = G_k(\iota_\mathcal{U}\oplus\iota_\mathcal{H})\psi(Z)u.
\end{eqnarray*}
Since the linear span of vectors of the form $\psi(Z)u,\ Z\in\Pi^{n_1\times n_1}\times\cdots\times\Pi^{n_N\times n_N},\ u\in\mathcal{U}$, is dense in $\mathcal{U\oplus H}$, the operator $G_k$ is $(\iota_\mathcal{U}\oplus\iota_\mathcal{H},\iota_{n_k}\otimes\iota_{\mathcal{M}_k})$-real, as desired.

(ii)$\Rightarrow$(i) Let $f$ satisfy (ii). Then by Theorem~\ref{thm:b-decomp} $f\in\mathcal{B}_{n_1,\ldots,n_N}(\mathcal{U})$. Let us show that the operator-valued linear function $A(Z)$ from \eqref{lp} is $\iota_{\mathcal U}\oplus\iota_{\mathcal H}$-real. Since $G_k$ is $(\iota_{\mathcal  U}\oplus\iota_{\mathcal H},\iota_{n_k}\otimes\iota_{\mathcal{M}_k})$-real, one has $G_k(\iota_{\mathcal U}\oplus\iota_{\mathcal H})=(\iota_{n_k}\otimes\iota_{\mathcal{M}_k})G_k$ and $(\iota_{\mathcal U}\oplus\iota_{\mathcal H})G_k^*=G_k^*(\iota_{n_k}\otimes\iota_{\mathcal{M}_k}),\ k=1,\ldots,N$. The latter equality follows from the fact that for every $h\in\mathbb{C}^{n_k}\otimes\mathcal{M}_k,\ x\in\mathcal{U\oplus H}$:
\begin{eqnarray*}
\left\langle(\iota_{\mathcal U}\oplus\iota_{\mathcal H})G_k^*h,x \right\rangle_{\mathcal{U\oplus H}} &=& \left\langle(\iota_{\mathcal U}\oplus\iota_{\mathcal H})x,G_k^*h \right\rangle_{\mathcal{U\oplus H}}\\
=\left\langle G_k(\iota_{\mathcal U}\oplus\iota_{\mathcal H})x,h \right\rangle_{\mathbb{C}^{n_k}\otimes\mathcal{M}_k} &=& \left\langle (\iota_{n_k}\otimes\iota_{\mathcal{M}_k})G_kx,h \right\rangle_{\mathbb{C}^{n_k}\otimes\mathcal{M}_k}\\
=\left\langle (\iota_{n_k}\otimes\iota_{\mathcal{M}_k})h,G_kx \right\rangle_{\mathbb{C}^{n_k}\otimes\mathcal{M}_k} &=& \left\langle G_k^*(\iota_{n_k}\otimes\iota_{\mathcal{M}_k})h,x \right\rangle_{\mathcal{U\oplus H}}.
\end{eqnarray*}
Therefore,
\begin{eqnarray*}
& & (\iota_{\mathcal U}\oplus\iota_{\mathcal H})A(\bar{Z})(\iota_{\mathcal U}\oplus\iota_{\mathcal H}) = \sum\limits_{k=1}^N(\iota_{\mathcal U}\oplus\iota_{\mathcal H})G_k^*(\overline{Z_k}\otimes I_{\mathcal{M}_k})G_k(\iota_{\mathcal U}\oplus\iota_{\mathcal H})\\ &=& \sum\limits_{k=1}^NG_k^*(\iota_{n_k}\otimes\iota_{\mathcal{M}_k})(\overline{Z_k}\otimes I_{\mathcal{M}_k})(\iota_{n_k}\otimes\iota_{\mathcal{M}_k})G_k = \sum\limits_{k=1}^NG_k^*(\iota_{n_k}\overline{Z_k}\iota_{n_k}\otimes I_{\mathcal{M}_k})G_k\\ &=& \sum\limits_{k=1}^NG_k^*(Z_k\otimes I_{\mathcal{M}_k})G_k
=A(Z).
\end{eqnarray*}
The latter is equivalent to the identities
\begin{eqnarray*}
\iota_\mathcal{U}a(\bar{Z})\iota_\mathcal{U}=a(Z), & \iota_\mathcal{U}b(\bar{Z})\iota_\mathcal{H}=b(Z),\\
\iota_\mathcal{H}c(\bar{Z})\iota_\mathcal{U}=c(Z), & \iota_\mathcal{H}d(\bar{Z})\iota_\mathcal{H}=d(Z).
\end{eqnarray*}
Since $\iota_\mathcal{H}^2=I_\mathcal{H}$ and
$$(\iota_\mathcal{H}d(\bar{Z})^{-1}\iota_\mathcal{H})\cdot(\iota_\mathcal{H}d(\bar{Z})\iota_\mathcal{H})=(\iota_\mathcal{H}d(\bar{Z})\iota_\mathcal{H})\cdot(\iota_\mathcal{H}d(\bar{Z})^{-1}\iota_\mathcal{H})=I_\mathcal{H},$$
one has
$$\iota_\mathcal{H}d(\bar{Z})^{-1}\iota_\mathcal{H}=(\iota_\mathcal{H}d(\bar{Z})\iota_\mathcal{H})^{-1}=d(Z)^{-1}.$$
Therefore,
\begin{eqnarray*}
f^\sharp(Z) &=& \iota_\mathcal{U}f(\bar{Z})\iota_\mathcal{U}=\iota_\mathcal{U}(a(\bar{Z})-b(\bar{Z})d(\bar{Z})^{-1}c(\bar{Z}))\iota_\mathcal{U}\\
&=& \iota_\mathcal{U}a(\bar{Z})\iota_\mathcal{U}-(\iota_\mathcal{U}b(\bar{Z})\iota_\mathcal{H})\cdot(\iota_\mathcal{H}d(\bar{Z})^{-1}\iota_\mathcal{H})\cdot(\iota_\mathcal{H}c(\bar{Z})\iota_\mathcal{U})\\
&=& a(Z)-b(Z)d(Z)^{-1}c(Z)=f(Z),
\end{eqnarray*}
and $f$ is $\iota_{\mathcal U}$-real. Thus, (i) follows from (ii).

(iv)$\Rightarrow$(i) Let (iv) hold. Then the operator $U=U^*=U^{-1}$ from \eqref{cal-u} is $((\iota_{n_1+\cdots +n_N}\otimes\iota_\mathcal{X})\oplus\iota_\mathcal{U})$-real, i.e.,
$$\left[\begin{array}{cc}
\iota_{n_1+\cdots +n_N}\otimes\iota_\mathcal{X} & 0\\
0 & \iota_\mathcal{U}\end{array}\right]\cdot\left[\begin{array}{cc}
A & B\\
C & D\end{array}\right]\cdot
\left[\begin{array}{cc}
\iota_{n_1+\cdots +n_N}\otimes\iota_\mathcal{X} & 0\\
0 & \iota_\mathcal{U}\end{array}\right]=\left[\begin{array}{cc}
A & B\\
C & D\end{array}\right].$$
This is equivalent to the following identities:
\begin{eqnarray*}
(\iota_{n_1+\cdots +n_N}\otimes\iota_\mathcal{X})A(\iota_{n_1+\cdots +n_N}\otimes\iota_\mathcal{X})=A, & (\iota_{n_1+\cdots +n_N}\otimes\iota_\mathcal{X})B\iota_\mathcal{U}=B,\\
 \iota_\mathcal{U}C(\iota_{n_1+\cdots +n_N}\otimes\iota_\mathcal{X})=C, & \iota_\mathcal{U}D\iota_\mathcal{U}=D.
 \end{eqnarray*}
 Moreover, for $W\in\mathbb{D}^{n_1\times n_1}\times\cdots\times\mathbb{D}^{n_N\times n_N}$ one has
\begin{eqnarray*}
\lefteqn{(\iota_{n_1+\cdots +n_N}\otimes\iota_\mathcal{X})(I_{\mathbb{C}^{n_1+\cdots +n_N}\otimes\mathcal{X}}-A(P(W)\otimes I_\mathcal{X}))(\iota_{n_1+\cdots +n_N}\otimes\iota_\mathcal{X})}\\
&=& (\iota_{n_1+\cdots +n_N}\otimes\iota_\mathcal{X})^2-(\iota_{n_1+\cdots +n_N}\otimes\iota_\mathcal{X})A(\iota_{n_1+\cdots +n_N}\otimes\iota_\mathcal{X})\\
&\times & (\iota_{n_1+\cdots +n_N}\otimes\iota_\mathcal{X})(P(W)\otimes I_\mathcal{X})(\iota_{n_1+\cdots +n_N}\otimes\iota_\mathcal{X})\\
&=& I_{\mathbb{C}^{n_1+\cdots +n_N}\otimes\mathcal{X}}-A(P(\overline{W})\otimes I_\mathcal{X}).
\end{eqnarray*}
Therefore,
\begin{eqnarray*}
(\iota_{n_1+\cdots +n_N}\otimes\iota_\mathcal{X})(I_{\mathbb{C}^{n_1+\cdots +n_N}\otimes\mathcal{X}}-A(P(W)\otimes I_\mathcal{X}))^{-1}(\iota_{n_1+\cdots +n_N}\otimes\iota_\mathcal{X})\\
=(I_{\mathbb{C}^{n_1+\cdots +n_N}\otimes\mathcal{X}}-A(P(\overline{W})\otimes I_\mathcal{X}))^{-1}
\end{eqnarray*}
(we already used an analogous argument above). Thus,
\begin{eqnarray*}
\mathcal{F}^\sharp(W) &=& \iota_\mathcal{U}\mathcal{F}(\overline{W})\iota_\mathcal{U}\\
&=& \iota_\mathcal{U}\left[D+C(P(\overline{W})\otimes I_\mathcal{X})(I_{\mathbb{C}^{n_1+\cdots +n_N}\otimes\mathcal{X}}-A(P(\overline{W})\otimes I_\mathcal{X}))^{-1}B\right]\iota_\mathcal{U}\\
&=& D+C(P(W)\otimes I_\mathcal{X})(I_{\mathbb{C}^{n_1+\cdots +n_N}\otimes\mathcal{X}}-A(P(W)\otimes I_\mathcal{X}))^{-1}B\\
&=& \mathcal{F}(W),\quad W\in\mathbb{D}^{n_1\times n_1}\times\cdots\times\mathbb{D}^{n_N\times n_N},
\end{eqnarray*}
i.e., $\mathcal{F}$ is $\iota_\mathcal{U}$-real. Applying the inverse double Cayley transform to $\mathcal{F}$, one can see that $f$ is $\iota_\mathcal{U}$-real on $\Pi^{n_1\times n_1}\times\cdots\times\Pi^{n_N\times n_N}$, and hence on $\Omega_{n_1,\ldots,n_N}$. Thus, (i) follows from (iv).

(iii)$\Rightarrow$(iv) Let $f$ satisfy \eqref{b-decomp} with holomorphic $(\iota_\mathcal{U},\iota_{n_k}\otimes\iota_{\mathcal{M}_k})$-real $L(\mathcal{U},\mathbb{C}^{n_k}\otimes\mathcal{M}_k)$-valued functions $\varphi_k$ on $\Pi^{n_1\times n_1}\times\cdots\times\Pi^{n_N\times n_N}$. As in the proof of Theorem~\ref{thm:image}, we get for $\mathcal{F}=\mathcal{C}(f)$ consecutively: identities \eqref{bb-id} and \eqref{bb-id-2} with holomorphic $L(\mathcal{U},\mathbb{C}^{n_k}\otimes\mathcal{M}_k)$-valued functions $\theta_k$ on $\Pi^{n_1\times n_1}\times\cdots\times\Pi^{n_N\times n_N}$ given by \eqref{theta-circ}, \eqref{theta} which are, moreover, $(\iota_\mathcal{U},\iota_{n_k}\otimes\iota_{\mathcal{M}_k})$-real; then identities \eqref{bb-id-1'} and \eqref{bb-id-2'}
with holomorphic $L(\mathbb{C}^{n_k}\otimes\mathcal{M}_k,\mathcal{U})$-valued functions $\widetilde{\theta_k}$ on $\Pi^{n_1\times n_1}\times\cdots\times\Pi^{n_N\times n_N}$  which are, moreover, $(\iota_{n_k}\otimes\iota_{\mathcal{M}_k},\iota_\mathcal{U})$-real. Without loss of generality, we consider all of $\mathcal{M}_k$'s equal, i.e., set $\mathcal{M}:=\bigoplus_{k=1}^N\mathcal{M}_k$ and regard
$$H_k^R=\theta_k\in L(\mathcal{U},\mathbb{C}^{n_k}\otimes\mathcal{M}),\quad H_k^L=\widetilde{\theta_k}\in L(\mathbb{C}^{n_k}\otimes\mathcal{M},\mathcal{U}),\quad \iota_\mathcal{M}:=\bigoplus_{k=1}^N\iota_{\mathcal{M}_k}.$$
Then $H_k^R$ is $(\iota_\mathcal{U},\iota_{n_k}\otimes\iota_\mathcal{M})$-real, and $H_k^L$ is $(\iota_{n_k}\otimes\iota_{\mathcal{M}},\iota_\mathcal{U})$-real, $k=1,\ldots,N$.

Set $\mathcal{X}:=\mathcal{M}$. Let us observe that the subspace $\mathcal{D}_0$ (and hence, ${\rm clos}(\mathcal{D}_0)$) is $((\iota_{n_1+\cdots +n_N}\otimes\iota_\mathcal{X})\oplus\iota_\mathcal{U})$-invariant. Indeed, for
$$x=\left[\begin{array}{c}
(P(W)\otimes I_\mathcal{X})H^R(W)\\
I_\mathcal{U}\end{array}\right]u+\left[\begin{array}{c}
H^R(W')\\
\mathcal{F}(W')\end{array}\right]u'\in\mathcal{D}_0,$$
 with some $u,u'\in\mathcal{U}$ and $W,W'\in\mathbb{D}^{n_1\times n_1}\times\cdots\times\mathbb{D}^{n_N\times n_N}$, we have
 \begin{eqnarray*}
 \lefteqn{((\iota_{n_1+\cdots +n_N}\otimes\iota_\mathcal{X})\oplus\iota_\mathcal{U})x}\\
 &=& \left[\begin{array}{c}
(\iota_{n_1+\cdots +n_N}P(W)\otimes \iota_\mathcal{X})H^R(W)\\
\iota_\mathcal{U}\end{array}\right]u+\left[\begin{array}{c}
(\iota_{n_1+\cdots +n_N}\otimes\iota_\mathcal{X})H^R(W')\\
\iota_\mathcal{U}\mathcal{F}(W')\end{array}\right]u'\\
 &=& \left[\begin{array}{c}
(P(\overline{W})\iota_{n_1+\cdots +n_N}\otimes \iota_\mathcal{X})H^R(W)\\
\iota_\mathcal{U}\end{array}\right]u+\left[\begin{array}{c}
H^R(\overline{W'})\iota_\mathcal{U}\\
\mathcal{F}(\overline{W}')\iota_\mathcal{U}\end{array}\right]u'\\
 &=& \left[\begin{array}{c}
(P(\overline{W})\otimes I_\mathcal{X})H^R(\overline{W})\\
I_\mathcal{U} \end{array}\right]\iota_\mathcal{U}u+\left[\begin{array}{c}
H^R(\overline{W'})\\
\mathcal{F}(\overline{W}')\end{array}\right]\iota_\mathcal{U}u\in\mathcal{D}_0.
  \end{eqnarray*}
Therefore, the subspace $$\mathcal{D}_0^\perp:=((\mathbb{C}^{n_1+\cdots +n_N}\otimes\mathcal{X})\oplus\mathcal{U})\ominus {\rm clos}(\mathcal{D}_0)$$ is also $((\iota_{n_1+\cdots +n_N}\otimes\iota_\mathcal{X})\oplus\iota_\mathcal{U})$-invariant. Indeed, for any $h_1\in\mathcal{D}_0,\ h_2\in\mathcal{D}_0^\perp$:
$$\left\langle ((\iota_{n_1+\cdots +n_N}\otimes\iota_\mathcal{X})\oplus\iota_\mathcal{U})h_2,h_1\right\rangle=\left\langle ((\iota_{n_1+\cdots +n_N}\otimes\iota_\mathcal{X})\oplus\iota_\mathcal{U})h_1,h_2\right\rangle=0,$$
thus $h_2 \in\mathcal{D}_0^\perp$ implies $((\iota_{n_1+\cdots +n_N}\otimes\iota_\mathcal{X})\oplus\iota_\mathcal{U})h_2\in\mathcal{D}_0^\perp$.

Now, it is easy to check that $U_0$ and therefore $\widetilde{U_0}$ are  $((\iota_{n_1+\cdots +n_N}\otimes\iota_\mathcal{X})\oplus\iota_\mathcal{U})$-real. Since $U=\widetilde{U_0}\oplus I_{\mathcal{D}_0^\perp}$, and $\mathcal{D}_0^\perp$ is $((\iota_{n_1+\cdots +n_N}\otimes\iota_\mathcal{X})\oplus\iota_\mathcal{U})$-invariant, the operator $U=U^*=U^{-1}$ is $((\iota_{n_1+\cdots +n_N}\otimes\iota_\mathcal{X})\oplus\iota_\mathcal{U})$-real, as required.

The proof is complete.
\end{proof}

\end{document}